\documentclass[11pt,intlimits]{amsart}

\usepackage{amsthm,amsmath,amssymb,amsopn,amsfonts,amscd,amsbsy}
\usepackage{mathtools}
\usepackage{array,caption}
\usepackage{url}
\usepackage{fullpage}

\theoremstyle{thmstyleone}
\newtheorem{theorem}{Theorem}

\newtheorem{proposition}[theorem]{Proposition}
\newtheorem{context}[theorem]{Context}
\newtheorem{corollary}[theorem]{Corollary}
\newtheorem{lemma}[theorem]{Lemma}

\theoremstyle{thmstyletwo}

\newtheorem{remark}{Remark}

\theoremstyle{thmstylethree}
\newtheorem{definition}{Definition}

\makeatletter

\@namedef{subjclassname@2020}{%
	\textup{2020} Mathematics Subject Classification}
\makeatother

\newcommand{\F}{\mathbb{F}}

\newcommand{\Z}{\mathbb{Z}}

\newcommand{\inert}{\text{inert}}
\newcommand{\ram}{\text{ram}}
\newcommand{\sm}{\text{sm}}
\newcommand{\spl}{\text{split}}
\newcommand{\tors}{\text{tors}}
\newcommand{\unr}{\text{unr}}

\newcommand{\calE}{\mathcal{E}}
\newcommand{\calO}{\mathcal{O}}
\newcommand{\calT}{\mathcal{T}}

\newcommand{\Aut}{\mathrm{Aut}}
\newcommand{\Gal}{\mathrm{Gal}}
\newcommand{\Spec}{\mathrm{Spec}}
\newcommand{\Tw}{\mathbf{Tw}}

\DeclareMathOperator{\ord}{ord}

\begin{document}

\title[$L$-Functions of Elliptic Curves Modulo Integers]{$L$-Functions of Elliptic Curves Modulo Integers}
\author{F\'elix Baril Boudreau}
\address{Department of Mathematics and Computer Science, University of Lethbridge, 4401 University Drive West, Lethbridge, T1K 3M4, Alberta, Canada}
\email{felix.barilboudreau@uleth.ca}

\begin{abstract}
	In 1985, Schoof devised an algorithm to compute zeta functions of elliptic curves over finite fields by directly computing the numerators of these rational functions modulo sufficiently many primes (see \cite{schoof_1985}). If $E/K$ is an elliptic curve with nonconstant $j$-invariant defined over a function field $K$ of characteristic $p \geq 5$, we know that its $L$-function $L(T,E/K)$ is a polynomial in $\Z[T]$ (see \cite[p.11]{katz_2002}). Inspired by Schoof, we study the reduction of $L(T,E/K)$ modulo integers. We obtain three main results. Firstly, if $E/K$ has non-trivial $K$-rational $N$-torsion for some integer $N$ coprime with $p$, we extend a formula for $L(T,E/K) \bmod N$ due to Hall (see \cite[p.133, Theorem 4]{hall_2006}) to all quadratic twists $E_f/K$ with $f \in K^\times \smallsetminus K^{\times 2}$. Secondly, without any condition on the $2$-torsion subgroup of $E(K)$, we give a formula for the quotient modulo $2$ of $L$-functions of any two quadratic twists of $E/K$. Thirdly, we use these results to compute the global root numbers of an infinite family of quadratic twists of an elliptic curve and in most cases find the exact analytic rank of each of these twists. We also illustrate that in favourable situations our second main result allows one to compute much more efficiently $L(T,E_f/K) \bmod 2$ than an algorithm of Baig and Hall (see \cite{baig_hall_2012}). Finally, we use our formulas to compute directly some degree $2$ $L$-functions.
\end{abstract}

\keywords{Elliptic Curve, Function Field, $L$-Function, Torsion Subgroup}

\subjclass[2020]{11G40}

\maketitle
\section{Introduction}\label{sec1}

\subsection{Motivation}
Let $k$ be a finite field whose cardinality $q$ is a power of a prime $p$ and let $E/k$ be an elliptic curve defined over $k$. The zeta function of $E/k$ is a rational function of the form
$$
Z(T,E/k) = \frac{L(T,E/k)}{(1-T)(1-qT)}.
$$
Here, $L(T,E/k)$ is the quadratic polynomial $1-\alpha T + qT^2$, with $\alpha$ the integer $1+q - \#E(k)$ and where $\#E(k)$ is the cardinality of the set of $k$-rational points $E(k)$. For $p \geq 5$, the elliptic curve $E/k$ is given by a Weierstrass equation of the form $y^2 = x^3 + ax + b$ with $a, b \in k$ satisfying $-16(4a^3 + 27b^2) \neq 0$. If we want to compute the coefficient $\alpha$ in practice, we can test all pairs $(u,v)$ in $k \times k$ and determine which ones satisfy the equation $v^2=u^3+au+b$. As the cardinality of $q$ grows, this method becomes impractical. In 1985, Schoof \cite{schoof_1985} side-stepped this point-counting method. He took advantage of the fact that $\alpha$ was bounded in absolute value by $2\sqrt{q}$ (see \cite[p.138, V.Theorem 1.1]{silverman_2009}) and devised an algorithm to directly compute the reduction of $\alpha$ modulo any prime $\ell$ different from $p$. If one performs this algorithm for sufficiently many primes $\ell_i$ so that their product $\prod_i \ell_i$ is strictly larger than $4\sqrt{q}$, then $\alpha$ is completely determined by its reduction modulo $\prod_i \ell_i$, thanks to the Chinese remainder theorem.

Suppose that $p \geq 5$ and let $C/k$ be a proper, smooth, and geometrically connected curve, with function field $K$. Let $E/K$ be an elliptic curve whose $j$-invariant is nonconstant. It follows from the work of Grothendieck, Deligne and others that the $L$-function of $E/K$ (defined in subsection \ref{LFunctionEllipticCurve}) is a polynomial in $1 + T \cdot \Z[T]$ (see \cite[p.11]{katz_2002}). In this paper, by analogy with Schoof, we begin a study of the problem of directly computing the reduction of $L(T,E/K)$ modulo an integer $N$ not divisible by $p$.
\subsection{Related Works}
The problem of explicitly computing the $L$-function of an elliptic curve defined over a function field is difficult and currently very few strategies to compute it exist. To our knowledge, they all take advantage of the observation that an elliptic curve $E/K$ can be seen as an elliptic surface $\calE/k$ over the finite field $k$ (see \cite[pp.250-252, Lecture 3]{ulmer_2011}). 

One can compute the $L$-function using point-counting methods. The only implemented point-counting method that we are aware of can be found in the computer algebra system Magma \cite{magma_1997} and is currently only available for the function field $K = k(t)$. The bulk of the work of the algorithm is to count points fiber-by-fiber on elliptic curves defined over finite fields that are finite extensions of $k$, up to some large enough degree. As the cardinality of $k$ grows, it becomes computationally difficult to compute $L$-functions with this technique.

In the direction of $p$-adic methods, one could consider the following idea, suggested by Andrew Sutherland (private communication) that we summarize here. In \cite{costa_harvey_kedlaya_2019}, Costa, Harvey and Kedlaya obtained a practical algorithm to compute zeta functions of varieties of some moderate dimension over finite fields which is much more efficient than naive-point counting. The idea would be to adapt their algorithm to compute the zeta function $Z(T,S/k)$ of the elliptic surface $S/k$ corresponding to the elliptic curve $E/K$. Let $Z(T,C/k)$ be the zeta function of the curve $C/k$ corresponding to the function field $K$ (see Definition \ref{DefinitionZetaFunction}). Sutherland's idea can be justified by the fact that the $L$-function $L(T,E/K)$ is (see \cite[p.257, Exercise 6.2]{ulmer_2011}), up to a finite product of rational functions indexed by the places of bad reduction, the quotient $Z(T,C/k)Z(qT,C/k)/Z(T,S/k)$. There are various algorithms to compute the zeta function of the curve $C/k$ (for example \cite{costa_harvey_kedlaya_2019}) and so to compute the $L$-function of $E/K$ it would essentially suffice to compute the zeta function of the elliptic surface. Costa pointed (private communication) that their method already applies if the elliptic surface is a toric $3$-fold, for example a surface in a weighted projective space.

Finally, since the $L$-function of an elliptic curve $E/K$ with nonconstant $j$-invariant is a polynomial with integer coefficients, one can explore another approach in spirit of Schoof's algorithm. The starting point of this exploration is the following key result by Hall \cite[p.133, Theorem 4]{hall_2006}.
\begin{theorem}\label{HallThm}
	Let $N \geq 2$ be an integer coprime with $q$ and let $\mathcal{T}$ be a finite subgroup of $E(K)$ of cardinality $N$. If $M^\text{sp}$, resp. $M^\text{ns}$, resp. $A$, denote the sets of places of $K$ over which $E$ has split multiplicative, resp. non-split multiplicative, resp. additive reduction, then
	\begin{align*}
		L(T,E/K) &\equiv Z(T,C/k)Z(qT,C/k) \times \prod\limits_{v \in M^\text{sp}} (1-q^{d_v}T^{d_v})\\
		&\times \prod\limits_{v \in M^\text{ns}}\frac{(1-T^{d_v})(1-q^{d_v}T^{d_v})}{(1+T^{d_v})}\\
		&\times \prod\limits_{v \in A}(1-T^{d_v})(1-q^{d_v}T^{d_v}) \bmod N.
	\end{align*}
\end{theorem}
It is worth mentioning that this formula was not used to compute $L$-functions in Hall's paper (see \cite[p.139, Theorem 12]{hall_2006} and \cite[p.142, Theorem 16]{hall_2006} for his main results). However, we illustrate in Proposition \ref{HallDegree2}, how Theorem \ref{HallThm} can be used to compute some degree $2$ $L$-functions.
\subsection{Overview of the Results}
Let $k$ be a finite field whose cardinality $q$ is a power of a prime $p \geq 5$ and let $C/k$ be a proper, smooth and geometrically connected curve with function field $K$. Let $E$ be an elliptic curve over $K$ with nonconstant $j$-invariant and for each $f \in K^\times \smallsetminus K^{\times 2}$ let $E_f/K$ be its quadratic twist (see Subsection \ref{QuadraticTwists}). Denote by $U$, $M$, $M^{\mathrm{sp}}, M^{\mathrm{ns}}$ and $A$ the subsets of closed points of $C$ corresponding to the loci of good, resp. multiplicative, resp. split multiplicative, resp. nonsplit multiplicative, resp. additive reduction of $E/K$. Let $f \in K^\times \smallsetminus K^{\times 2}$ and let $K_f$ be the quadratic extension of $K$ generated by a chosen square root of $f$. We decorate by a subscript $f$ the corresponding sets $U_f$, $M_f$, $M_f^{\mathrm{sp}}, M_f^{\mathrm{ns}}$ and $A_f$ when the reduction types of the places of $K$ are considered with respect to the quadratic twist $E_f/K$ of $E/K$. If $v$ is a place of $K$ (equivalently a closed point of $C$), then we write $k_v$ for the residue field of $v$ and $d_v$ for the degree of the finite extension $k_v/k$. Let $M_{\unr}^{\mathrm{sp}}$ (resp. $M_{\unr}^{\mathrm{ns}}$) be the subset of $M^{\mathrm{sp}}$ (resp. of $M^{\mathrm{ns}}$) whose elements are unramified in $K_f$. Finally, let $M_{\spl}$ and $M_{\inert}$ be the subsets of places of $M$ which are respectively split and inert in $K_f$ and let $A_\unr, A_\ram^{\text{gd}}$ be subsets of places of $A$ which are respectively unramified in $K_f$ and ramified in $K_f$ over which $E/K$ has potentially good reduction.
Our first main result, Theorem \ref{ThmA}, gives, under the conditions of Theorem \ref{HallThm} a formula for $L(T,E_f/K) \bmod N$ for any of the infinitely many $f \in K^\times \smallsetminus K^{\times 2}$.
\begin{theorem}[Theorem \ref{Thm1}]\label{ThmA} Suppose that $E(K)$ contains a subgroup $\mathcal{T}$ of order $N \geq 2$ with $N$ coprime with $q$. Let $\chi$ be the nontrivial character associated to the quadratic extension $K_f/K$ and let $L(T,\chi)$ be its corresponding Artin $L$-function (see Definition \ref{ArtinLFunction}). If $p \geq 5$, then
	\begin{align*}
		L(T,E_f/K) \equiv &L(T,\chi)L(qT,\chi) R(T) \bmod N,
	\end{align*}
	where
	$$
	R(T) =
	\prod_{v \in M^{\text{sp}}_{\text{unr}}} (1 + \alpha_v q^{d_v}T^{d_v}) \prod_{v \in M^{\text{ns}}_{\text{unr}}} \frac{(1+\alpha_v T^{d_v})(1+\alpha_v q^{d_v}T^{d_v})}{1-\alpha_vT^{d_v}},
	$$
	with
	$$
	\alpha_v =
	\begin{cases}
		1 & \text{ if } v \in M_{\text{inert}},\\
		-1 & \text{ if } v \in M_{\text{split}}.\\
	\end{cases}
	$$
\end{theorem}
\begin{remark}
	In Theorem \ref{Thm1} we also provide formulas for $L(T,E_f/K) \bmod N$ when $N \in \{2,3,4\}$. As a consequence, Corollary \ref{Cor1} gives a formula for $L(T,E_f/K) \bmod 2$ that is slightly more precise than the one in \cite[p.133, Theorem 4]{hall_2006}.
\end{remark}
\begin{remark}
	It is important to point that even if $E(K)$ has a subgroup of order $N$, the group $E_f(K)$ does not need to have a subgroup of order $N$. Indeed, one can show that for $N$ odd, the direct sum of $N$-torsion subgroups $E(K)[N] \oplus E_f(K)[N]$ of $E(K)$ and $E_f(K)$ is isomorphic to $N$-torsion subgroup $E(K_f)[N]$ of $E(K_f)$. Moreover, in \cite[p.133, Theorem 4]{hall_2006} only gives $L(T,E/K) \bmod N$, while Theorem \ref{ThmA} gives $L(T,E_f/K) \bmod N$ any of the infinitely many $f \in K^\times \smallsetminus K^{\times 2}$.
\end{remark}
The next result conveniently provides $\deg L(T,E_f/K)$, the degree of the $L$-function of a quadratic twist $E_f/K$, and its global root number $\varepsilon(E_f/K)$, in terms of the places of bad reduction of $E/K$, the decomposition types of these places in $K_f$, and the genus $g(C)$ of $C$.
\begin{proposition}[Proposition \ref{DegreeAndRootNumberQuadraticTwists}]\label{PropA}
	Let $k'$ be the field of constants of $K_f$.
	\begin{enumerate}
		\item[(i)] The degree of $L(T,E_f/K)$ equals
		$$
		\deg(M_\unr) + \deg(A_\unr) \text{ if } [k':k] = 2,
		$$
		$$
		2\left( 2g(C) - 2 \right) + 2 \hspace{-0.5cm}\sum\limits_{\substack{v \in \mid C\mid  \\ \ord_v(f) \equiv 1 \bmod 2}} \hspace{-0.8cm}d_v + \deg(M_\unr) + \deg(A_\unr) - \deg(A_\ram^\text{gd}) \text{ if } [k':k] = 1.
		$$
		\item[(ii)] If $E(K)$ has a subgroup of order $N$ with $N \geq 2$ an integer coprime with $q$, then the global root number $\varepsilon(E_f/K)$ satisfies
		$$
		\varepsilon(E_f/K) \equiv
		\begin{cases}
			(-1)^{\# M_\unr^{\text{sp}}+\# M_\inert^{\text{ns}}} q^{-\deg(A_\unr)} \bmod N & \text{ if } [k':k] = 2,\\
			(-1)^{\# M_\unr^{\text{sp}}+\# M_\inert^{\text{ns}}} q^{-\deg(A_\unr) + \deg(A_\ram^{\text{gd}})} \bmod N & \text{ if } [k':k] = 1.
		\end{cases}
		$$
	\end{enumerate}
\end{proposition}
In our second main result we provide, unconditionally on the $2$-torsion subgroup of $E(K)$, the ratio modulo $2$ of the $L$-functions of any two quadratic twists of $E/K$. This ratio is well-defined. Indeed, let $f_1,f_2 \in K^\times$ and let $E_{f_1}$ and $E_{f_2}$ be the corresponding quadratic twists of $E$. Since $E$ has nonconstant $j$-invariant $j(E)$, then $L(T,E/K) \in 1 + T \cdot \Z[T]$ (see \cite[p.11]{katz_2002}). By Lemma \ref{DiscriminantAndJInvariantQuadraticTwist}, we have $j(E) = j(E_{f_i})$ and so $L(T,E_{f_i}/K) \in 1 + T \cdot \Z[T]$, for $i = 1,2$. Therefore, $L(T,E_{f_1}/K)/L(T,E_{f_2}/K) \in 1 + T \cdot \Z[[T]]$ and $L(T,E_{f_1}/K)/L(T,E_{f_2}/K) \bmod 2$ is a well-defined element of $1 + T \cdot (\Z/2\Z)[[T]]$.
\begin{theorem}[Theorem \ref{Thm2}] \label{ThmB}
	Let $f_1,f_2 \in K^\times$ and let 
	$$
	U_{f_1,f_2} := \left(\mid U_{f_1}\mid  \cap \mid U_{f_2}\mid \right) \cup \left(M_{f_1} \cap M_{f_2}\right) \cup \left(A_{f_1} \cap A_{f_2}\right).$$
	Then,
	$$
	\frac{L(T,E_{f_1}/K)}{L(T,E_{f_2}/K)} \equiv \prod_{v \not \in U_{f_1,f_2}} \frac{L_v(T^{d_v},E_{f_2}/K)}{L_v(T^{d_v},E_{f_1}/K)} \bmod{2}.
	$$
\end{theorem}
In particular, when $K = k(t)$, we have the following.
\begin{corollary}[Corollary \ref{Cor2}]\label{CorA}
	Let $q$ be a power of a prime $p \geq 5$ and let $K = k(t)$. Suppose that $E/K$ has semistable reduction everywhere, with multiplicative reduction at $t=\infty$. Let $\Delta \in k[t]$ be the discriminant of a minimal Weierstrass equation for $E$ and let $f \in k[t]$ be a square-free polynomial of even degree which is coprime with $\Delta$. We have
	\begin{equation}\label{UnconditionalModulo2}
		\frac{L(T,E_f/K)}{L(T,E/K)} \equiv \prod_{v \mid f} (1 + a_vT^{d_v} + T^{2d_v}) \bmod 2,
	\end{equation}
	where $a_v := 1 + q^{d_v} - \#E_v(k_v)$. Here, $E_v$ is the elliptic curve over $k_v$ obtained from $E$ by reducing modulo $v$ the coefficients of that minimal Weierstrass equation, and $\#E_v(k_v)$ is the cardinality of its set of $k_v$-rational points. See Corollary \ref{Cor2} for a more general expression.
\end{corollary}
A highlight of Equation (\ref{UnconditionalModulo2}) is its computability as we illustrate in Subsection \ref{Complexity}. Given an elliptic curve $E/K$ with $L(T,E/K) = 1$, we discuss the efficiency to compute the reduction $L(T,E_f/K)$ modulo $2$, when $f$ is a monic irreducible polynomial of even degree $d_f$, which is coprime to $\Delta$. More precisely, we show that employing Equation (\ref{UnconditionalModulo2}) is by far superior to computing $L(T,E/K) \bmod 2$ than using a suitable version of the algorithm of Baig and Hall in \cite{baig_hall_2012}. The reason is roughly the following. Since $char(k) \neq 2,3$, we can write a Weierstrass equation for $E/K$ in the form $y^2 = x^3 + Ax + B$ for some $A, B \in k[t]$. The comparison of the two approaches to compute $L(T,E_f/K) \bmod 2$ amounts to decide whether the image of $x^3 + Ax + B$ in some finite fields is irreducible or not. With Corollary \ref{Cor2}, we have $L(T,E_f/K) \equiv 1 + a_fT^{d_f} + T^{2d_f} \bmod 2$ and we then only need to compute the reduction of $x^3 + Ax + B$ in the finite extension $k[t]/(f)$ of $k$ of degree $d_f$. In general, the algorithm of Baig and Hall requires to test, for most monic irreducible polynomials $v$ in $k[t]$ of degree $d_v \leq 2 d_f$, the irreducibility of the image of $x^3 + Ax + B$ in the extensions of finite fields of $k$ defined by these polynomials $v$. As $d_f$ grows, the number of finite fields to consider for the algorithm of Baig and Hall becomes considerable (see Subsection \ref{Complexity} for more details).

In Subsection \ref{ComputationLFunction}, we discuss an analogue approach to Schoof's algorithm to compute polynomial $L$-functions of elliptic curves with nonconstant $j$-invariants. Although this analogous approach is currently in a primitive state, Theorem \ref{HallThm}, Theorem \ref{ThmA} and Theorem \ref{Thm2} can be seen as first concrete steps. We also make the following elementary observation.
\begin{proposition}[Proposition \ref{OddAnalyticRank} and Corollary \ref{FirstCasesOddAnalyticRank}]\label{PropB}
	If $E/K$ is an elliptic curve with nonconstant $j$-invariant, with odd analytic rank, with a polynomial $L$-function of degree $d$ and coefficients $a_0 = 1, a_1, \ldots, a_n \in \Z$, then
	\begin{equation}\label{EquationOddAnalyticRank}
		L(T,E/K) = (1-qT)\left( 1 + \sum_{n=1}^{d-2}\left( \sum_{m=0}^n a_m q^{n-m} \right)T^n + q^{d-1}T^{d-1} \right). 
	\end{equation}
\end{proposition}
As an illustration, we conclude Subsection \ref{ComputationLFunction} by using Theorem \ref{HallThm}, Theorem \ref{ThmA}, Proposition \ref{PropA} and Corollary \ref{CorA} to compute some degree $2$ $L$-functions of elliptic curves.
\begin{proposition}[Proposition \ref{HallDegree2}]
	Let $k$ be the field with $5$ elements and let $f$ be a rational function of the form $f \in \{ P(t), 1/Q(t), P(t)/Q(t), \text{ with } P(t), Q(t) \in k[t] \text{ coprime and of degree } 1 \}$. An elliptic curve $E/k(t)$ given by  a Weierstrass equation
	$$
	E/k(t) : y^2 = (1-a)xy - by = x^3 - bx^2,
	$$
	with $a = (10-2f)/(f^2-9)$ and $b = -2(f-1)^2(f-5)/(f^2-9)^2$, has torsion subgroup $E(k(t))_\tors$ isomorphic to $\Z/2\Z \times \Z/6\Z$. Moreover, its $L$-function equals
	$$
	L(T,E/k(t)) = (1-5T)(1+5T).
	$$
\end{proposition}
\begin{proposition}[Proposition \ref{ComputationQuadraticTwists}]
	Suppose that $k$ is the field with $7$ elements and let $f \in k[t]$ with degree $1$. An elliptic curve $E/k(t)$ given by a Weierstrass equation
	$$
	E/k(t): y^2 + xy + \frac{f^2+f+1}{3(f+2)^3} = x^3
	$$
	has torsion subgroup $E(k(t))_\tors$ isomorphic to $(\Z/3\Z)^2$, bad reduction $M^{\text{sp}}$ and $L$-function equal to $1$. Moreover, for any finite place $h \in M^{\text{sp}}$, the $L$-function of the quadratic twist $E_h/k(t)$ is
	$$
	L(T,E_h/k(t)) = (1-7T)(1+7T).
	$$
\end{proposition}
Finally, in Subsection \ref{StudyInfiniteFamily}, we study, over an arbitrary field $k$ of characteristic $p \geq 5$, an infinite family of quadratic twists of the so-called universal elliptic curve $E$ over the function field $K$ of the modular curve $X_1(3)/k$. We obtain our third main theorem (detailed in Theorem \ref{Thm3}) by using Theorem \ref{ThmA}, Proposition \ref{PropA} and Corollary \ref{CorA} to compute the reduction modulo $2$ and modulo $3$ of the $L$-functions of the quadratic twists of that family, compute all the global root numbers and either give the exact analytic rank of the quadratic twists or a nontrivial upper bound on it.
\begin{theorem}[Theorem \ref{Thm3}]\label{ThmC}
	Let $k$ be a finite field of cardinality $q$ such that $q \equiv 5, 11 \bmod 12$. Consider the elliptic curve $E/k(t)$ given by the Weierstrass equation
	\begin{equation}\label{IntroWeierstrassEquationX1(3)}
		y^2 + 3xy + (1-t^3)y = x^3
	\end{equation}
	and let $\Delta$ be the discriminant of this equation. The point $(0,0) \in E(k(t))$ has order $3$. If $f \in k[t]$ is a monic square-free polynomial of degree $d_f$, which is coprime to $\Delta$, then $L(T,E_f/K)$ has even degree. For $a \in \{0,1\}$, we write $f(a) \in k^2$? to ask whether $f(a)$ is a square (Yes) or not (No) in $k$. We have

	\begin{table}[h]
		\centering
		\begin{tabular}{|c|c|c|}
			\hline
			$f(0) \in k^2$? & $f(1) \in k^2$? & $L(T,E_f/K) \bmod 3$\\
			\hline
			Yes & Yes & $L(T,K_f)L(-T,K_f)(1-T)(1+T)(1+T^2)$\\
			\hline
			Yes & No & $L(T,K_f)L(-T,K_f)(1-T)^2(1+T^2)$\\
			\hline
			No & Yes & $L(T,K_f)L(-T,K_f)(1-T)(1+T)^3$\\
			\hline
			No & No & $L(T,K_f)L(-T,K_f)(1-T)(1+T)(1+T^2)$\\
			\hline
		\end{tabular}
		\caption{Reduction $L(T,E_f/K) \bmod 3$}
\label{tb1Mod3}
\end{table}

	The global root number $\varepsilon(E_f/K)$ is given in Table \ref{tb2RootandRank} below.
	
	\begin{table}[h]
		\centering
		\begin{tabular}{ |c|c|c|c|}
			\hline
			$f(0) \in k^2$? & $f(1) \in k^2$? & $\varepsilon(E_f/K)$ & $\mathrm{rank}^{an}(E_f/K)$\\
			\hline
			Yes & Yes & $1$ & $0$ \\
			\hline
			Yes & No & $1$ & $0$ \\
			\hline
			No & Yes & $-1$ & $\leq 3$ \\
			\hline
			No & No & $-1$ & $1$ \\
			\hline
		\end{tabular}
		\caption{Global Root Number and Analytic Rank Using Reduction Modulo $3$}
		\label{tb2RootandRank}
	\end{table}

	Moreover, let $\mathrm{Jac}(C_f)(k)$ be the set of $k$-rational points of the Jacobian of the (normalization of the) curve $C_f/k$ associated to the field $K_f$. If the $3$-torsion subgroup of $\mathrm{Jac}(C_f)(k)$ is trivial, then the fourth column of Table \ref{tb2RootandRank} gives the analytic rank, $\mathrm{rank}^{an}(E_f/K)$ of the elliptic curve $E_f/K$ or an upper bound on it.
	
	If $f$ is irreducible of odd degree, then we have the following.
	\begin{table}[h]
	\centering
		\begin{tabular}{ |c|c|c|}
			\hline
			$a_\infty \bmod 2$ & $a_f \bmod 2$ & $\mathrm{rank}^{an}(E_f/K)$ \\
			\hline
			$0$ & $0$ & $\leq 4$\\
			\hline
			$0$ & $1$ & $\leq 2$\\
			\hline
			$1$ & $0$ & $\leq 2$\\
			\hline
			$1$ & $1$ & $0$\\
			\hline
		\end{tabular}
		\caption{Analytic Rank using Reduction Modulo $2$}
\label{tb3Mod2}
\end{table}
\end{theorem}

This paper is organized as follows. In Section \ref{NotationPrerequisites}, we introduce the notation used throughout the paper and recall some key concepts about zeta functions of curves and about quadratic twists and $L$-functions of elliptic curves. In Section \ref{ArbitraryGenus} we prove our first and second main results, Theorem \ref{Thm1} and Theorem \ref{Thm2}, as well as Corollary \ref{Cor1} and Proposition \ref{DegreeAndRootNumberQuadraticTwists}. In Section \ref{FunctionFieldofGenus0}, we specialize to $K = k(t)$ the statements of Theorem \ref{Thm1} and Theorem \ref{Thm2} in Proposition \ref{SimplificationsThm1} and Corollary \ref{Cor2}, respectively. In Section \ref{Application}, we discuss three applications of the results of Section \ref{ArbitraryGenus} and Section \ref{FunctionFieldofGenus0}. More precisely, we discuss in Subsection \ref{Complexity} the efficiency of Corollary \ref{Cor2} to compute $L(T,E_f/K) \bmod 2$ and compare it with an algorithm of Baig and Hall. In Subsection \ref{ComputationLFunction}, we discuss an analogue approach of Schoof's algorithm to compute polynomial $L$-functions, observe an elementary identity on $L$-functions of elliptic curves with odd analytic rank (see Proposition \ref{OddAnalyticRank}) and compute some degree $2$ $L$-functions with this Schoof's like approach. In Subsection \ref{StudyInfiniteFamily}, we prove our third main result, Theorem \ref{Thm3}. Finally, in Section \ref{Appendix}, we prove some auxiliary results mentioned in Section \ref{NotationPrerequisites} and in Section \ref{ArbitraryGenus} for which we do not know a reference in the literature.
\section{Notation and Prerequisites}\label{NotationPrerequisites}
In this section, we set up the notation that we use throughout this paper, introduce some basic concepts and prove some results which are certainly well-known but for which we could not find a proof in the literature.

Let $k$ be a finite field of cardinality a power $q$ of a prime $p \geq 5$. Let $C/k$ be a smooth, proper and geometrically connected curve $C$ defined over $k$ and of genus $g := g(C)$. Let $\mid C\mid $ be the the set of closed points of $C/k$ and let $K$ be its function field, which also has genus $g = g(K)$. The field $K$ is a finite algebraic extension of $k(t)$ for some transcendental element $t \in K$ over $k$. Given a field $L$, we fix an algebraic closure $\overline{L}$.

Now suppose that $K$ is the function field of a curve $C/k$ as introduced above. If $v$ is a place of $K$, then $K_v$ denotes the completion of $K$ at $v$ and we let $\ord_v$ be the discrete valuation associated to $v$. Let $\calO_v$ be the discrete valuation ring $ \calO_v = \{ \alpha \in K_v : \ord_v(\alpha) \geq 0 \}$. The ring $\calO_v$ is a local ring and we choose for generator of its unique maximal ideal $\mathfrak{m}_v = \{ \alpha \in K_v : \ord_v(\alpha) > 0 \}$ a uniformizer $\pi_v$, normalized so that $\ord_v(\pi_v) = 1$. We denote by $k_v$ the associated residue field $\calO_v/\mathfrak{m}_v$, which is a finite extension of degree $d_v$ of $k$. Given an elliptic curve $E/K$, we let $E_v/k_v$ be the reduced curve defined over $k_v$ and we let $\Phi_v$ by the local group scheme of components defined over $k_v$. If $X$ is a scheme and $L$ is a field, then $X(L)$ denotes the set of $L$-rational points of $X$. If $G$ is a group and a positive integer $N$, then $G[N]$ denotes the $N$-torsion subgroup of $G$. Finally, if $S$ is a finite set, then $\# S$ denotes its cardinality.
\begin{definition}\label{DefinitionZetaFunction}
	Let $C/k$ be a proper, smooth, and geometrically connected curve of genus $g$. We define the zeta function of $C/k$ as the Euler product
	$$
	Z(T,C/k) := \prod_{v \in \mid C\mid } (1-T^{d_v})^{-1}.
	$$
	This zeta function can also be expressed as rational function
	$$
	Z(T,C/k) = \frac{L(T,C/k)}{(1-T)(1-qT)},
	$$
	where $L(T,C/k)$ is a polynomial in $1 + T \cdot \Z[T]$ of degree $2g$ (see \cite[p.53, Theorem 5.9]{rosen_2002}).
	
\end{definition}
\begin{definition}
	Let $E/K$ be an elliptic curve. Consider inside the locus of bad reduction $Z \subset C$ all closed points $v$ according to their reduction types with respect to $E/K$. We write $M$, resp. $M^{\text{sp}}$, resp. $M^{\text{ns}}$, resp. $A$, for the set of points of multiplicative reduction, resp. of split multiplicative reduction, resp. of nonsplit multiplicative reduction, resp. of additive reduction. If $S$ is any of these above sets, we sometimes abuse notation and consider $S$ as the divisor $\sum_{v \in S}d_v\cdot v$ and in this case define $\deg(S) = \sum_{v \in S}d_v$.
\end{definition}
The following result is certainly known, but we could not find it in the literature. Compare with the proof of \cite[p.132, Lemma 3]{hall_2006}.
\begin{lemma}\label{PrimeTopTorsion}
	Let $K_v$ be a field complete with respect to a non-Archimedean place $v$ whose residue field $k_v$ has characteristic $p \geq 0$. If $E/K_v$ is an elliptic curve with additive reduction over $v$, then, the prime-to-$p$-torsion subgroup of $E(K_v)$ injects into the group $\Phi_v(k_v)$.
\end{lemma}
\begin{proof}
	See Lemma \ref{AppendixPrimeTopTorsion} in the Appendix, Section \ref{Appendix}.
\end{proof}

\subsection{Quadratic Twists}\label{QuadraticTwists}
We now briefly recall the notion of quadratic twists of elliptic curves and state some of their properties. For more details, we refer the reader to \cite[X.2 and X.5]{silverman_2009}. Let $(E/K,O_E)$ be an elliptic curve $E$ defined over $K$ together with a rational point $O_E \in E(K)$.
\begin{definition}
	A \textit{twist of the elliptic curve} $E/K$ is the data of a second elliptic curve $E'/K$ together with an isomorphism of elliptic curves $\phi: E'/\overline{K} \to E/\overline{K}$ (so that $\phi$ sends $O_{E'}$ to $O_E$). On the set $\calT(E/K)$ of twists of the elliptic curve $E/K$, we introduce the equivalence relation $\sim_K$ that two twists are equivalent if they are $K$-isomorphic and denote by $\Tw(E/K)$ the resulting quotient set.
\end{definition}
The set $\Aut(E)$ of $\overline{K}$-automorphisms of the group structure of the elliptic curve $E/K$ forms an abelian group (see \cite[p.69]{silverman_2009}). One obtains a canonical set bijection (see \cite[p.319, X.2.2(c) and X.2.3]{silverman_2009})
from $\Tw(E/K)$ to the cohomology group $\mathrm{H}^1(\Gal(\overline{K}/K),\Aut(E))$.
Now, let
$$
n := \begin{cases}
	2 & \text{ if } j(E) \neq 0, 1728,\\
	4 & \text{ if } j(E) = 1728,\\
	6 & \text{ if } j(E) = 0.
\end{cases}
$$
Then, there is a canonical isomorphism of $\Gal(\overline{K}/K)$-modules $\mu_n \xrightarrow{\simeq} \Aut(E)$ (see \cite[p.104, III.10.2]{silverman_2009}). As a consequence of Hilbert's Theorem 90, we have a canonical group isomorphism $K^\times/K^{\times n} \xrightarrow{\simeq} \mathrm{H}^1(\Gal(\overline{K}/K),\mu_n)$ (see \cite[p.420, Appendix B, Proposition 2.5 (c)]{silverman_2009}). Therefore, the set $\Tw(E/K)$ is canonically isomorphic to the set $K^\times/K^{\times n}$. In this paper, an elliptic curve $E/K$ is assumed to have nonconstant $j$-invariant and so we are only concerned with the case $n=2$. In this situation, if we choose a Weierstrass equation of the form $y^2 = x^3 + ax + b$ for $E/K$ and an element $f \in K^\times$, then by \cite[p.343, X.5.4]{silverman_2009}, a representative $E_f/K$ of a class in $\Tw(E/K)$ corresponding to $f \bmod K^{\times 2}$ has a Weierstrass equation of the form $y^2 = x^3 + f^2ax + f^3 b$. We call the elliptic curve $E_f/K$ the \textit{quadratic twist of $E/K$ by $f$}. When $f = 1$, the curve $E_f/K$ is precisely the ``original'' elliptic curve $E/K$.
An immediate application of the definitions of discriminant of a Weierstrass equation and of $j$-invariant of an elliptic curve yields the following result.
\begin{lemma}\label{DiscriminantAndJInvariantQuadraticTwist}
	Let $L$ be any field of characteristic different from $2$ and $3$.
	Let $E/L$ be an elliptic curve, let $f \in L^\times$ and let $E_f/L$ be the corresponding quadratic twist. Choose a Weierstrass equation of the form $W: y^2 = x^3 + ax + b$ for $E/L$ and let $\Delta_W$ be its discriminant. The corresponding Weierstrass equation for $E_f/L$ is $ W_f : y^2 = x^3 + f^2 ax + f^3b$ and we denote by $\Delta_{W_f}$ its discriminant. We have $\Delta_{W_f} = f^6\Delta_W$ and in particular $j(E_f) = j(E)$.
\end{lemma}
\subsection{$L$-Function of an Elliptic Curve}\label{RecallsLFunctions}
Finally, we define the $L$-function of an elliptic curve defined over a function field.
Let $C$ be a smooth, proper and geometrically connected curve defined over a finite field $k$ of characteristic $p \geq 5$ of genus $g$. Let $K$ be the function field of $C/k$ and let $E/K$ be an elliptic curve. For each place $v$ of $K$, we define the following integer.
\begin{equation}\label{TraceOfFrobenius}
	a_v :=
	\begin{cases}
		1 + q^{d_v} - \#E_v(k_v) & \text{ if } E/K \text{ has good reduction over }v,\\ 
		-1 & \text{ if } E/K \text{ has split multiplicative reduction over }v,\\
		1 &\text{ if } E/K \text{ has nonsplit multiplicative reduction over }v,\\
		0 &\text{ if } E/K \text{ has additive reduction over }v.
	\end{cases}
\end{equation}
\begin{definition}\label{ConductorsEllipticCurve}
	For each place $v$ of $K$, we define, following \cite[p.381, 10.2(b)]{silverman_1994} the \textit{local exponent of the conductor} at $v$ to be the integer
	$$
	n_v = \begin{cases}
		0 & \text{ if } E/K \text{ has good reduction over } v,\\
		1 & \text{ if } E/K \text{ has multiplicative reduction over } v,\\
		2 & \text{ if } E/K \text{ has additive reduction over } v.\\
	\end{cases}
	$$
	The \textit{global conductor} of $E/K$ is the divisor $N(E/K) := \sum_{v \in \mid C\mid } n_v[v]$. Its degree equals
	\begin{equation}\label{DegreeLFunction}
		\deg N(E/K) = \sum_v n_v \deg(v) = \deg(M) + 2\deg(A),
	\end{equation}
	where $M$ and $A$ respectively denote the divisors of places of $K$ over which $E/K$ has multiplicative and additive reduction. 
\end{definition}
\begin{definition}\label{LFunctionEllipticCurve}
	Given a formal indeterminate $T$, set
	\begin{equation}\label{EulerFactorLFunction}
		L_v(T^{d_v},E/K) :=
		\begin{cases}
			1 - a_vT^{d_v} + q^{d_v}T^{2d_v} & \text{ if } E/K \text{ has good reduction over } v,\\
			1-a_vT^{d_v} & \text{ if } E/K \text{ has bad reduction over } v.\\
		\end{cases}
	\end{equation}
	We define $L$-function of $E/K$ is defined to be the Euler product
	\begin{equation}\label{EulerProduct}
		L(T,E/K) := \prod_{v \in \mid C\mid } L_v(T^{d_v},E/K)^{-1}.
	\end{equation}
\end{definition}
Now assume that the $j$-invariant of $E/K$ is nonconstant. It follows from the work of Grothendieck, Deligne and others (see \cite[p.11]{katz_2002}) that $L(T,E/K)$ is a polynomial in $1 + T \cdot \Z[T]$ and has degree \cite[p.232, Theorem 9.3]{ulmer_2011}
\begin{equation}\label{DegreeLFunction}
	d = 2(2g - 2) + \deg N(E/K).
\end{equation}
Moreover, the leading coefficient of $L(T,E/K)$ is $\varepsilon(E/K)q^d$, where $\varepsilon(E/K)$ is the global root number of $E/K$. For this paper, it suffices to know that $\varepsilon(E/K) = 1$ precisely when $\mathrm{rank}^{an}(E/K)$, the analytic rank of $E/K$, i.e., the order of vanishing of $L(T,E/K)$ at $T = q^{-1}$, is even and otherwise $\varepsilon(E/K) = -1$.
\section{Function Fields of Arbitrary Genus}\label{ArbitraryGenus}
We assume throughout this section the following context.
\begin{context}\label{ContextArbitraryGenus}
	Let $C/k$ be a proper, smooth and geometrically connected curve of genus $g$ defined over a finite field $k$ of characteristic $p \geq 5$ and let $K$ be its function field. Let $f \in K^\times \smallsetminus K^{\times 2}$ be a non-square element and let $E/K$ be an elliptic curve with nonconstant $j$-invariant. Modulo $K^{\times 2}$, the element $f$ determines a unique quadratic extension $K_f$ of $K$ as well as a unique quadratic twist $E_f/K$. Let $\theta \in \overline{K}$ be such that $\theta^2 = f$ and so that $K_f = K(\theta)$. Given a place $v$ of $K$, we let $f_v$ be the image of $f$ via a chosen field injection $K \xhookrightarrow{} K_v$ and we let $\theta_v$ be the image of $\theta$ via the corresponding field injection $\overline{K} \xhookrightarrow{} \overline{K_v}$. This yields a $K$-embedding $\tau: K_f \xhookrightarrow[]{} \overline{K_v}$. By abuse of notation, let $v$ be the canonical extension to $K_v$ of the place $v$ and let $\overline{v}$ be the unique extension of the place $v$ of $K_v$ to $\overline{K_v}$. Then we have an extension $w = \overline{v} \circ \tau$ of $v$ to $K_f$. Since the map $\tau$ is continuous, it extends uniquely to a continuous $K$-embedding $\tau: K_{f,w} \xhookrightarrow[]{} \overline{K_v}$ and since $K_f/K$ is finite, then $K_{f,w}$ is the completion of $K_f$ at $w$. One then sees that $K_{f,w}$ is the compositum field $K_fK_v$ \cite[p.161, II.8]{neukirch_1999} and the latter equals $K_v(\theta_v)$ by definition of compositum of fields.
\end{context}

Let $\mid C\mid _\text{split}$ (resp. $\mid C\mid _\text{inert}$, resp. $\mid C\mid _\text{ram}$) be the subset of places of $K$ which are split (resp. inert, resp. ramified) in $K_f$, so that $\mid C\mid  = \mid C\mid _\text{split} \cup \mid C\mid _\text{inert} \cup \mid C\mid _\text{ram}$ is a disjoint union. Define $\mid C\mid _\text{unr}:=\mid C\mid _\text{split} \cup \mid C\mid _\text{inert}$ and if $S \in \{M^{\text{sp}},M^{\text{ns}},A\}$ and $r \in \{\text{inert, ram, split, unr}\}$, set $S_r:=S\cap\mid C\mid _r$.

Also, write $A_{\text{ram}}^{\text{sp}}$, respectively $A_{\text{ram}}^{\text{ns}}$, respectively $A_{\text{ram}}^{\text{gd}}$, for the places of $K$ which are ramified in $K_f$ and such that, with respect to the elliptic curve $E/K$, are places of additive reduction which are more precisely potentially multiplicative split, respectively potentially multiplicative nonsplit, respectively potentially good. Finally, let $A_{\text{ram}}^m$ be the union $A_{\text{ram}}^{\text{sp}} \cup A_{\text{ram}}^{\text{ns}}$ and let $A_{\text{ram}}^o$ be the complement in $A_{\text{ram}}$ of $A_{\text{ram}}^m \cup A_{\text{ram}}^{\text{gd}}$.

For $f \in K^\times \smallsetminus K^{\times 2}$ we decorate the corresponding sets as $S_{r,f}$, $A_{f,\text{ram}}^{\text{sp}}$, $A_{f,\text{ram}}^{\text{ns}}$, $A_{f,\text{ram}}^{\text{gd}}$, $A_{f,\text{ram}}^m$ and $A_{f,\text{ram}}^o$, when these sets refer to places whose reduction types are considered with respect to the quadratic twist $E_f/K$.
\begin{remark}\label{RamificationValuation}
	We have $\mid C\mid _\text{ram} = \{ v \in \mid C\mid  : \ord_v(f) \equiv 1 \bmod 2\}$.
\end{remark}
\subsection{General Results}
The following lemma holds in general and we use it in Subsection \ref{StudyInfiniteFamily}.
\begin{lemma}\label{NumeratorZetaFunctionAt-1}
	We have
	$$
	L(-1,C/k) = q^gL(-q^{-1}, C/k).
	$$
	In particular, if $N$ is a positive integer such that $q \equiv -1 \bmod N$, then
	$$
	L(-1,C/k) \equiv (-1)^g \#\mathrm{Jac}(C)(k) \bmod N,
	$$
	where $\#\mathrm{Jac}(C)(k)$ is the number of $k$-rational points of the Jacobian of $C/k$.
\end{lemma}
\begin{proof}
	This follows directly from the functional equation $L(T,C/k) = q^g T^{2g}L( q^{-1}T^{-1}, C/k)$ (see \cite[p.54, Proof of Theorem 5.9]{rosen_2002}).
\end{proof}

The following result gives a formula for the degree of $L(T,E_f/K)$ and also extend \cite[p.134, Corollary 5]{hall_2006} to quadratic twists. Moreover, the degree of the $L$-function of $E_f/K$ and the global root number of $E_f/K$ are expressed in terms of the data of $E/K$ and of the decomposition behaviour of the places of $K$ in $K_f$.

\begin{proposition}\label{DegreeAndRootNumberQuadraticTwists}
	Let $k'$ be the field of constants of the quadratic extension $K_f$ of $K$.
	\begin{enumerate}
		\item[(i)] The degree of $L(T,E_f/K)$ equals
		$$
		\deg(M_\unr) + \deg(A_\unr) \text{ if } [k':k] = 2,
		$$
		$$
		2\left( 2g(C) - 2 \right) + 2 \hspace{-0.5cm}\sum\limits_{\substack{v \in \mid C\mid  \\ \ord_v(f) \equiv 1 \bmod 2}} \hspace{-0.8cm}d_v + \deg(M_\unr) + \deg(A_\unr) - \deg(A_\ram^\text{gd}) \text{ if } [k':k] = 1.
		$$
		\item[(ii)] If $E(K)$ has a subgroup of order $N$ with $N \geq 2$ an integer coprime with $q$, then the global root number $\varepsilon(E_f/K)$ satisfies
		$$
		\varepsilon(E_f/K) \equiv
		\begin{cases}
			(-1)^{\# M_\unr^{\text{sp}}+\# M_\inert^{\text{ns}}} q^{-\deg(A_\unr)} \bmod N & \text{ if } [k':k] = 2,\\
			(-1)^{\# M_\unr^{\text{sp}}+\# M_\inert^{\text{ns}}} q^{-\deg(A_\unr) + \deg(A_\ram^{\text{gd}})} \bmod N & \text{ if } [k':k] = 1.
		\end{cases}
		$$
	\end{enumerate}
\end{proposition}
\begin{proof}
	If $[k':k] = 2$, then $K_f$ is the compositum $k' K$, the extension $K_f/K$ is unramified and $g(K_f) = g(K) = g(C)$ (see \cite[p.114, Theorem 3.6.3 (a),(b)]{stichtenoth_2009}).
	
	Now assume that $k' = k$. Since $p \geq 5$, then the extension $K_f/K$ is tamely ramified. It follows from Remark \ref{RamificationValuation} and from the Riemann-Hurwitz formula that
	$$
	g(K_f) = 2g(K) - 1 + \frac{1}{2}  \sum_{\substack{v \in \mid C\mid  \\ \ord_v(f) \equiv 1 \bmod 2}} d_v.
	$$
	Let $U^f$, $M^f$ and $A^f$ be the divisors of places of $K_f$ over which the elliptic curve $E$ has respectively good, multiplicative and additive reduction. For any quadratic extension $K_f/K$, we have $L(T,E_f/K) = L(T,E/K) L(T,E_f/K)$.
		\begin{enumerate}
		\item[(i)] From Equation (\ref{DegreeLFunction}) it follows that
		\begin{align*}
			\deg(L(T,E_f/K)) &= \deg(L(T,E/K_f)) - \deg(L(T,E/K))\\
			&= 4\left( g(K_f) - g(K) \right) + \left( \deg(M^f) - \deg(M) \right)\\
			&+ 2 \left( \deg(A^f) - \deg(A) \right).
		\end{align*}
		Observe that if $v \in M_\text{inert} \cup A_\text{inert}$ has degree $d_v$, then the degree of the place $w \in M^f$ equals $2d_v$. Moreover, if $v \in A_\ram^{\text{gd}}$, then a place $w$ of $K_f$ lying over $v$ belongs to $\mid U^f \mid$. Finally, recall that $M_\text{unr} = M_\text{split} \cup M_{\text{inert}}$.  Then,
				\begin{align*}
			\deg(M^f) &= 
			\begin{cases}
				2\deg(M_\text{unr}) & \text{ if } [k':k] = 2,\\
				2\deg(M_\text{unr}) + \deg(M_\text{ram})& \text{ if } [k':k] = 1,\\
			\end{cases}\\
			\deg(A^f) &= 
			\begin{cases}
				2\deg(A_\text{unr}) & \text{ if } [k':k] = 2,\\
				2\deg(A_\text{unr}) + \deg(A^m_\text{ram}) + \deg(A^o_\text{ram})& \text{ if } [k':k] = 1.\\
			\end{cases}\\
		\end{align*}
	Thus,
	\begin{align*}
		\deg(M^f) - \deg(M) &= \deg(M_\unr) \text{ for any } [k':k] \in \{1,2\}, \text{ and}\\
		\deg(A^f) - \deg(A) &=
		\begin{cases}
			\deg(A_\text{unr}) & \text{ if } [k':k] = 2,\\
			\deg(A_\text{unr}) - \deg(A_\ram^\text{gd}) & \text{ if } [k':k] = 1.
		\end{cases}
	\end{align*}
	Therefore, $\deg\left( L(T,E_f/K) \right)$ equals
	$$
	\deg(M_\unr) + \deg(A_\unr) \text{ if } [k':k] = 2,
	$$
	$$
	2\left( 2g(C) - 2 \right) + 2 \sum\limits_{\substack{v \in \mid C\mid  \\ \ord_v(f) \equiv 1 \bmod 2}} d_v + \deg(M_\unr) + \deg(A_\unr) - \deg(A_\ram^\text{gd}) \text{ if } [k':k] = 1.
	$$
	\item[(ii)] We have $\varepsilon(E/K_f) = \varepsilon(E/K)\varepsilon(E_f/K)$. Since $N$ divides $\# E(K)$, then $N$ also divides $E(K_f)$. Therefore, from \cite[p.134, Corollary 5]{hall_2006} we have
	$$
	\varepsilon(E/K_f) \equiv (-1)^{\# M^{f,\text{sp}}} q^{-\deg(A^f)} \bmod N \text{ and } \varepsilon(E/K) \equiv (-1)^{\# M^{\text{sp}}} q^{-\deg(A)} \bmod N.
	$$
	Moreover, observe that $\# M^{f,\text{sp}} = 2\# M^{\text{sp}} + \# M^{\text{ns}}_\inert$. Therefore,
	\begin{align*}
		\varepsilon(E_f/K) &\equiv (-1)^{\# M^{f,\text{sp}} -\# M^{\text{sp}}} q^{\deg(A) - \deg(A^f)} \bmod N\\
		&\equiv
		\begin{cases}
			(-1)^{\# M_\unr^{\text{sp}}+\# M_\inert^{\text{ns}}} q^{-\deg(A_\unr)} \bmod N & \text{ if } [k':k] = 2,\\
			(-1)^{\# M_\unr^{\text{sp}}+\# M_\inert^{\text{ns}}} q^{-\deg(A_\unr) + \deg(A_\ram^{\text{gd}})} \bmod N & \text{ if } [k':k] = 1.
		\end{cases}
	\end{align*}
		\end{enumerate}
\end{proof}
The following lemma is needed in our proofs of Theorem \ref{Thm1} and Theorem \ref{Thm2}. It is certainly known but we could not find it in the literature. We refer the reader to Lemma \ref{AppendixLemmaKodairaSymbols} in the Appendix, Section \ref{Appendix}, for our proof of the result.
\begin{lemma}\label{LemmaKodairaSymbols}
	The reduction types of $E/K$ and of $E_f/K$ are related as follows.
	\begin{enumerate}
		\item[(a)] If $v \in \mid C\mid _{\text{unr}}$,
		\begin{enumerate}
			\item[(a1)] then $v \in \mid U\mid  \cap \mid U_f\mid $ or $v \in M \cap M_f$ or $v \in A \cap A_f$.
			\item[(a2)] Moreover, if $v \in M$, then $v \in M^{\text{sp}} \cap M_f^{\text{sp}}$ or $M^{\text{ns}} \cap M_f^{\text{ns}}$ if and only if $v \in M_{\text{split}}$.
		\end{enumerate}  
		\item[(b)] If $v \in \mid C\mid _{\text{ram}}$, then
	\begin{table}[h]
	\centering
			\begin{tabular}{ |c|c|c|c|c|c| }
				\hline
				$E/K_v$ & $\mid U\mid _{\text{ram}}$ & $M_{\text{ram}}$ & $A_\text{ram}^{\text{gd}}$ & $A_\text{ram}^m$ & $A_\text{ram}^{o}$\\
				\hline
				$E_f/K_v$ & $A_{\text{ram}}^\text{gd}$  & $A_{f,\text{ram}}^m$ & $\mid U_f\mid _{\text{ram}}$ & $M_{f,\text{ram}}$ & $A_{f,\text{ram}}^{o}$\\
				\hline
			\end{tabular}.
		\caption{Reduction Types After Ramified Twist}
		\label{tb4RamifiedTwist}
		\end{table}
	\end{enumerate}
\end{lemma}
\begin{definition}\label{ArtinLFunction}
	Given $f \in K^\times \smallsetminus K^{\times 2}$, we let $\chi$ be the quadratic character associated to the quadratic extension $K_f/K$. For each place $v$ of $K$, set 
	$$
	L_v(T^{d_v},\chi) :=
	\begin{cases}
		1-T^{d_v} & \text{ if } v \in \mid C\mid _{\text{split}},\\
		1+T^{d_v} & \text{ if } v \in \mid C\mid _{\text{inert}},\\
		1                & \text{ if } v \in \mid C\mid _{\text{ram}}.\\
	\end{cases}
	$$
	The \textit{Artin} $L$-function of $\chi$ is the formal Euler product
	$$
	L(T,\chi) := \prod_{v \in \mid C\mid } L_v(T^{d_v},\chi)^{-1}.
	$$
\end{definition}
Suppose that $E(K)$ contains a subgroup of order $N \geq 2$ with $N$ coprime with $q$. Under this assumption, our first main result gives an explicit formula for $L(T,E_f/K)$ modulo $N$.
\begin{theorem}\label{Thm1}
	If $E(K)$ contains a subgroup $\mathcal{T}$ of order $N \geq 2$ with $N$ coprime with $q$, then,
	\begin{align*}
		L(T,E_f/K) \equiv &L(T,\chi)L(qT,\chi) \times R(T) \times \gamma_{\text{ram}} \bmod N,
	\end{align*}
	where
	$$
	R(T) \equiv
	\begin{cases}
		\prod_{v \in M^{\text{sp}}_{\text{unr}}} (1 + \alpha_v q^{d_v}T^{d_v}) \times \prod_{v \in M^{\text{ns}}_{\text{unr}}} \frac{(1+\alpha_v T^{d_v})(1+\alpha_v q^{d_v}T^{d_v})}{1-\alpha_vT^{d_v}} &\text{ if } N \geq 5,\\
		\prod_{v \in M^{\text{sp}}_{\text{unr}}} (1 + \alpha_v q^{d_v}T^{d_v}) \times \prod_{v \in M^{\text{ns}}_{\text{unr}}} \frac{(1+\alpha_v T^{d_v})(1+\alpha_v q^{d_v}T^{d_v})}{1-\alpha_vT^{d_v}} \\
		\times \prod\limits_{v \in A_{\text{unr}}} (1+\alpha_vT^{d_v})(1+\alpha_v q^{d_v}T^{d_v}) &\text{ if } N \in \{2,3,4\},
	\end{cases}
	$$
	with
	$$
	\alpha_v =
	\begin{cases}
		1 & \text{ if } v \in M_{\text{inert}} \cup A_{\text{inert}},\\
		-1 & \text{ if } v \in M_{\text{split}} \cup A_{\text{split}},\\
	\end{cases}
	$$
	and $\gamma_{\text{ram}} = 1$ if $N \geq 5$ and
	\begin{align*}
		\gamma_{\text{ram}} &= \prod\limits_{v \in A_{\text{ram}}^{\text{sp}}} (1 - T^{d_v})^{-1} \times \prod\limits_{v \in A_{\text{ram}}^{\text{ns}}} (1 + T^{d_v})^{-1}\\
		&\times \prod\limits_{v \in A_{\text{ram}}^{\text{gd}}} (1 - T^{d_v})^{-1}(1 + T^{d_v})^{-1}
	\end{align*}
	if $N \in \{2,3,4\}$.
	More precisely, for $N \in \{2,3,4\}$ the possible Kodaira symbols of additive reduction with respect to $E/K_v$ which contribute to $L(T,E_f/K) \bmod N$ are as follows.
\begin{table}[h]
\centering
		\begin{tabular}{ |c|c|c| }
			\hline
			$N$ & $A_\text{unr}$ & $A_{\text{ram}}$\\
			\hline
			$2$ &$III, III^*, I_n^*, n \geq 0$& $I_n^*, n \geq 0$\\
			\hline
			$3$ & $IV, IV^*$ & $I_n^*, n \geq 0$ \\
			\hline
			$4$ &$I_n^*, n \geq 0$ & $I_n^*, n \geq 0$\\
			\hline
		\end{tabular}
	\caption{Kodaira Symbols of Additive Reduction Contributing to $L(T,E_f/K) \bmod N$}
	\label{tb5PossibleKodairaSymbols}
\end{table}
\end{theorem}

Before proving Theorem \ref{Thm1}, we recall and prove the following useful lemma.

\begin{lemma}\label{LemAdditive}
	Let $K_v$ be a complete field with respect to some non-Archimedean valuation and with finite residue field of characteristic $p$. Let $E/K_v$ be an elliptic curve with additive reduction and let $\mathcal{T}$ be a torsion subgroup of $E(K_v)$ of order prime to $p$. Then $\mathcal{T}$ has order at most $4$. If it has order at least $2$, then we have the following classification.
	\begin{table}[h]
	\centering
		\begin{tabular}{ |c|c| }
			\hline
			Isomorphism Type of $\mathcal{T}$ & Possible Kodaira Symbols\\
			\hline
			$\Z/2\Z$ & $III, III^*, I_n^*$ with $n \geq 0$, $I^*_{m,2}$ with $m \geq 0$\\
			\hline
			$\Z/3\Z$ & $IV, IV^*$\\
			\hline
			$\Z/4\Z$ & $I_{2n+1}^*$, with $n \geq 1$\\
			\hline
			$(\Z/2\Z)^2$ & $I_{2n}^*$, with $n \geq 0$\\
			\hline
		\end{tabular}
	\caption{Local Torsion and Possible Additive Kodaira Symbols}
	\label{tb6LocalAdditiveReduction}
	\end{table}
\end{lemma}
\begin{remark}
	In the notation of \cite[p.497, 10.2.24]{liu_2002}, the symbol $I_{m,2}^*, m \geq 0$ is the Kodaira symbol $I_m^*, m \geq 0$ where the special fiber of the corresponding minimal model only has $2$ irreducible components of multiplicity $1$ defined over the residue field $k_v$.
\end{remark}
\begin{proof}(of Lemma \ref{LemAdditive})
	By Lemma \ref{PrimeTopTorsion}, the prime-to-$p$ torsion of $E(K_v)$ injects into $\Phi_v(k_v)$. The first part of the statement follows from the fact that $\calT$ is then identified with a subgroup $\Phi_v(k_v)$ and the latter has order at most $4$ (see \cite[p.362, IV.9.2(d)]{silverman_1994}). To show the second part of the statement, we obtain the following classification for the isomorphism type of the group $\Phi_v(k_v)$, based on the tables in \cite[p.497, 10.2.24]{liu_2002}.

	\begin{table}[h]
		\centering
		\begin{tabular}{ |c|c| }
			\hline
			Isomorphism Type of $\Phi_v(k_v)$ & Possible Kodaira Symbols\\
			\hline
			$\Z/2\Z$ & $III, III^*$, $I_{m,2}^*, m \geq 0$\\
			\hline
			$\Z/3\Z$ & $IV, IV^*$\\
			\hline
			$\Z/4\Z$ & $I_{2n+1}^*$ with $n \geq 1$\\
			\hline
			$(\Z/2\Z)^2$ & $I_{2n}^*$ with $n \geq 0$\\
			\hline
		\end{tabular}
	\caption{Possible Kodaira Symbols for a Local Component Group of Additive Reduction}
	\label{tb7LocalComponentGroupAndKodaira}
	\end{table}

	Finally, in order to include $I_n^*$ with $n \geq 0$, where $\Phi_v(k_v) = \Phi_v(\overline{k_v})$, as possible Kodaira symbols for which $\calT$ is isomorphic to $\Z/2\Z$, we observe that the groups $\Z/4\Z$ and $(\Z/2\Z)^2$ have subgroups of order $2$.
\end{proof}
\begin{proof}{ (of Theorem \ref{Thm1}) }
	The idea of the proof is as follows. We consider the quotient
	$$
	\frac{L(T,E_f/K)}{L(T,\chi)L(qT,\chi)} \equiv \prod_{v \in \mid C\mid } \frac{L_v(T^{d_v},\chi)L_v(q^{d_v}T^{d_v},\chi)}{L_v(T^{d_v},E_f/K)} \bmod N
	$$
	and analyze the ratios of Euler factors according to the subsets $\mid C\mid _{\text{split}}$, $\mid C\mid _{\text{inert}}$, $\mid C\mid _{\text{ram}}$ of $\mid C\mid $.
	
	We now prove the result. Using Lemma \ref{LemmaKodairaSymbols}, we obtain $L_v(T^{d_v},E_f/K)$ in terms of the reduction type of $E/K$ over $v \in \mid C\mid $.
	
	First suppose that $v \in \mid C\mid _{\text{split}}$. Then, we obtain the following table.
	
	\begin{table}[h]
		\centering
		\begin{tabular}{ |c|c|c| }
			\hline
			Reduction type of $E/K_v$ & $L_v(T^{d_v},E_f/K) \bmod N$ & $\frac{L_v(T^{d_v},\chi)L_v(q^{d_v}T^{d_v},\chi)}{L_v(T^{d_v},E_f/K)} \bmod N$\\
			\hline
			$\mid U\mid _{\text{split}}$ & $(1-T^{d_v})(1-q^{d_v}T^{d_v})$ & $1$\\
			\hline
			$M^{sp}_{\text{split}}$ & $1-T^{d_v}$ & $1-q^{d_v}T^{d_v}$\\
			\hline
			$M^{ns}_{\text{split}}$ & $1+T^{d_v}$ & $\frac{(1-T^{d_v})(1-q^{d_v}T^{d_v})}{1+T^{d_v}}$\\
			\hline
			$A_{\text{split}}$ & $1$ & $(1-T^{d_v})(1-q^{d_v}T^{d_v})$\\
			\hline
		\end{tabular}
	\caption{Ratio of Euler Factors When the Place Splits in $K_f$}
	\label{tb7RatioSplit}
	\end{table}

	Second suppose that $v \in \mid C\mid _{\text{inert}}$. Then, we obtain the following table.

	\begin{table}[h]
		\centering
		\begin{tabular}{ |c|c|c| }
			\hline
			Reduction type of $E/K_v$ & $L_v(T^{d_v},E_f/K) \bmod N$ & $\frac{L_v(T^{d_v},\chi)L_v(q^{d_v}T^{d_v},\chi)}{L_v(T^{d_v},E_f/K)} \bmod N$\\
			\hline
			$\mid U\mid _{\text{inert}}$ & $(1-T^{d_v})(1-q^{d_v}T^{d_v})$ & $1$\\
			\hline
			$M^{sp}_{\text{inert}}$ & $(1+T^{d_v})$ & $1+q^{d_v}T^{d_v}$\\
			\hline
			$M^{ns}_{\text{inert}}$ & $(1-T^{d_v})$ & $\frac{(1+T^{d_v})(1+q^{d_v}T^{d_v})}{1-T^{d_v}}$\\
			\hline
			$A_{\text{inert}}$ & $1$ & $(1+T^{d_v})(1+q^{d_v}T^{d_v})$\\
			\hline
		\end{tabular}
	\caption{Ratio of Euler Factors When the Place is Inert in $K_f$}
	\label{tb8RatioSplit}
	\end{table}

	Third suppose that $v \in \mid C\mid _{\text{ram}}$. In this case, $L_v(T^{d_v},\chi)L_v(q^{d_v}T^{d_v},\chi) = 1$. Then, we obtain the following table.

	\begin{table}[h]
	\centering
		\begin{tabular}{ |c|c|c| }
			\hline
			Reduction type of $E/K_v$ & $L_v(T^{d_v},E_f/K) \bmod N$ & $\frac{L_v(T^{d_v},\chi)L_v(q^{d_v}T^{d_v},\chi)}{L_v(T^{d_v},E_f/K)} \bmod N$\\
			\hline
			$\mid U\mid _{\text{ram}}$ & $1$ & $1$\\
			\hline
			$M_{\text{ram}}$ & $1$ & $1$\\
			\hline
			$A_{\text{ram}}^o$& $1$ & $1$\\
			\hline
			$A_{\text{ram}}^{\text{sp}}$& $(1-T^{d_v})^{-1}$ & $(1-T^{d_v})$\\
			\hline
			$A_{\text{ram}}^{\text{ns}}$& $(1+T^{d_v})^{-1}$ & $(1+T^{d_v})$\\
			\hline
			$A_{\text{ram}}^{\text{gd}}$& $(1-T^{d_v})(1+q^{d_v}T^{d_v})$ & $(1-T^{d_v})^{-1}(1+q^{d_v}T^{d_v})^{-1}$\\
			\hline
		\end{tabular}
	\caption{Ratio of Euler Factors When the Place Ramifies in $K_f$}
	\label{tb10RatioRamify}
	\end{table}

	Only the last line of Table \ref{tb10RatioRamify} requires a justification. Suppose that $E/K$ has additive potentially good reduction over $v$. Then, $E_f/K$ has good reduction over $v$ and by the theorem of semistable reduction, $E_f/K_f$ has good reduction over the place $w$ of $K_f$ lying above $v$. Since $E$ and $E_f$ are isomorphic over $K_f$, then $E_f(K_f)[N] = E(K_f)[N]$. Moreover, we have $E(K_f)[N] \supset E(K)[N] \neq \{ O \}$ by assumption and so $L_w(T^{d_w},E_f/K_f) \equiv (1-T^{d_w})(1-q^{d_w}T^{d_w}) \bmod N$ by \cite[p.132, Lemma 3]{hall_2006}. Since $v$ ramifies in $K_f$, then $d_v = d_w$. Therefore, $\#E_{f,v}(k_v) = \#E_w(k_w)$ and so $L_v(T^{d_v},E_f/K) \equiv (1-T^{d_v})(1-q^{d_v}T^{d_v}) \bmod N$.
	
	Now, $E(K)$ has a subgroup $\calT$ of order $N \geq 2$ with $N$ coprime to $q$. For each place $v$ of $K$, $E(K)$ injects into $E(K_v)$ and so does $\calT$. If $v$ is a place of additive reduction, then by Lemma \ref{LemAdditive}, $\calT$ has order at most $4$. In particular, for $N \geq 5$, the elliptic curve $E/K$ has no place of additive reduction. Therefore, for $N \geq 5$ we ignore the rows regarding places of additive reduction in the previous tables. 
\end{proof}

The following corollary of Theorem \ref{Thm1} is a more precise version of \cite[p.133, Theorem 4]{hall_2006} in the case $N = 2$. More precisely, we note that places of additive reduction for $E_f/K$ which ramify in $K_f$ do not contribute to the reduction $L(T,E_f/K) \bmod 2$.
\begin{corollary}[of Theorem \ref{Thm1}]\label{Cor1}
	If $E(K)$ has nontrivial $2$-torsion, then
	$$
	L(T,E_f/K) \equiv Z(T,C/k)^2 \times \prod\limits_{v \in M_f} (1-T^{d_v}) \times \prod_{v \in A_\text{f,unr}} (1-T^{d_v})^2 \bmod 2.\\
	$$
\end{corollary}
\begin{proof}
	Since $q$ is odd, it is congruent to $1$ modulo $2$. Moreover, $-1 \equiv 1 \bmod 2$ and so we do not distinguish between split and nonsplit multiplicative reduction modulo $2$.
	
	First, suppose that $v \in \mid C\mid _{\text{split}}$. Then, we obtain the following table.

	\begin{table}[h]
		\centering
		\begin{tabular}{ |c|c| }
			\hline
			Reduction type of $E/K_v$ & $\frac{L_v(T^{d_v},\chi)L_v(q^{d_v}T^{d_v},\chi)}{L_v(T^{d_v},E_f/K)} \bmod 2$\\
			\hline
			$\mid U\mid _{\text{split}}$ &  $1$\\
			\hline
			$M_{\text{split}}$ &  $1-T^{d_v}$\\
			\hline
			$A_{\text{split}}$& $(1-T^{d_v})^2$\\
			\hline
		\end{tabular}
	\caption{Ratio Modulo 2 When the Place Splits in $K_f$}
	\label{tb11Mod2RatioSplits}
    \end{table}

	Second, suppose that $v \in \mid C\mid _{\text{inert}}$. Then, we obtain the following table.

	\begin{table}[h]
		\centering
		\begin{tabular}{ |c|c| }
			\hline
			Reduction type of $E/K_v$ & $\frac{L_v(T^{d_v},\chi)L_v(q^{d_v}T^{d_v},\chi)}{L_v(T^{d_v},E_f/K)} \bmod 2$\\
			\hline
			$\mid U\mid _{\text{inert}}$ & $1$\\
			\hline
			$M_{\text{inert}}$ & $1-T^{d_v}$\\
			\hline
			$A_{\text{inert}}$ & $(1-T^{d_v})^2$\\
			\hline
		\end{tabular}
		\caption{Ratio Modulo 2 When the Place is Inert in $K_f$}
	\label{tb12Mod2RatioInrty}
	\end{table}

	Third, suppose that $v \in \mid C\mid _{\text{ram}}$. Then, we obtain the following table.

	\begin{table}[h]
		\centering
		\begin{tabular}{ |c|c|}
			\hline
			Reduction type of $E/K_v$ & $\frac{L_v(T^{d_v},\chi)L_v(q^{d_v}T^{d_v},\chi)}{L_v(T^{d_v},E_f/K)} \bmod 2$\\
			\hline
			$\mid U\mid _{\text{ram}}$ & $1$ \\
			\hline
			$M_{\text{ram}}$ & $1$ \\
			\hline
			$A_{\text{ram}}^o$& $1$ \\
			\hline
			$A_{\text{ram}}^m$& $(1-T^{d_v})^{-1}$\\
			\hline
			$A_{\text{ram}}^{\text{gd}}$& $(1-T^{d_v})^{-2}$\\
			\hline
		\end{tabular}
		\caption{Ratio Modulo 2 When the Place Ramifies in $K_f$}
	\label{tb13Mod2RatioRamify}
	\end{table}

	Therefore,
	\begin{align*}
		\frac{L(T,E_f/K)}{L(T,\chi)L(qT,\chi)} &\equiv \prod_{v \in M_{\text{unr}}} (1-T^{d_v}) \times \prod_{v \in A_{\text{unr}}} (1-T^{d_v})^2\\
		&\times \prod_{v \in A_{\text{ram}}^m} (1-T^{d_v})^{-1} \times \prod_{v \in A_{\text{ram}}^{\text{gd}}} (1-T^{d_v})^{-2} \bmod 2.
	\end{align*}
	Observe that $L(T,\chi) \equiv Z(T,C/k) \times \prod_{v \in \mid C\mid _\text{ram}} (1-T^{d_v}) \bmod 2$. Hence,
	\begin{align*}
		L(T,E_f/K) &\equiv Z(T,C/k)^2 \times \prod_{v \in A_{\text{ram}}^{\text{gd}}} (1-T^{d_v})\\
		&\times \prod_{v \in M_{\text{unr}}} (1-T^{d_v}) \times \prod_{v \in A_{\text{unr}}} (1-T^{d_v})^2\\
		&\equiv Z(T,C/k)^2 \times \prod\limits_{v \in M_f} (1-T^{d_v}) \times \prod_{v \in A_\text{f,unr}} (1-T^{d_v})^2 \bmod 2.
	\end{align*}
\end{proof}
The next lemma is used in the proof of Theorem \ref{Thm2}.
\begin{lemma}\label{LemmaLegendreSymbol}
	For each  $v \in \mid U\mid _\unr = \mid U\mid  \cap \mid U_f\mid $, we have $a_v + a_{f,v} \equiv 0 \bmod 2$, where $a_v$ and $a_{f,v}$ are defined in Equation (\ref{TraceOfFrobenius}).
\end{lemma}
\begin{proof}
	See Lemma \ref{LemmaLegendreSymbolAppendix}, in the Appendix, Section \ref{Appendix}.
\end{proof}
Let $f_1,f_2 \in K^\times$ and let $E_{f_1}/K$ and $E_{f_2}/K$ be the corresponding quadratic twists of $E/K$. Since $E/K$ has nonconstant $j$-invariant $j(E)$, then $L(T,E/K) \in 1 + T \cdot \mathbb{Z}[T]$ (see \cite[p.11]{katz_2002}). Since $j(E) = j(E_{f_i})$ for $i \in \{1,2\}$ by Lemma \ref{DiscriminantAndJInvariantQuadraticTwist}, then $L(T,E_{f_i}/K) \in 1 + T \cdot \mathbb{Z}[T]$ for $i \in \{1,2\}$. Therefore, $L(T,E_{f_1}/K)/L(T,E_{f_2}/K) \in 1 + T \cdot \mathbb{Z}[[T]]$ and the expression $L(T,E_{f_1}/K)/L(T,E_{f_2}/K) \bmod 2$ is a well-defined element of $1 + T \cdot (\Z/2\Z)[[T]]$. In the second main result of this paper, we write explicitly the reduction modulo $2$ of this quotient.
\begin{theorem}\label{Thm2}
	Let $f_1,f_2 \in K^\times$ and let 
	$$
	U_{f_1,f_2} := \left(\mid U_{f_1}\mid  \cap \mid U_{f_2}\mid \right) \cup \left(M_{f_1} \cap M_{f_2}\right) \cup \left(A_{f_1} \cap A_{f_2}\right).$$
	Then,
	$$
	\frac{L(T,E_{f_1}/K)}{L(T,E_{f_2}/K)} \equiv \prod_{v \not \in U_{f_1,f_2}} \frac{L_v(T^{d_v},E_{f_2}/K)}{L_v(T^{d_v},E_{f_1}/K)} \bmod{2}.
	$$
\end{theorem}
\begin{proof}
	The ratio of Euler products modulo $2$ is congruent to the product of the ratios modulo 2 of the Euler factors:
	$$
	\frac{L(T,E_{f_1}/K)}{L(T,E_{f_2}/K)} \equiv \prod_v \frac{L_v(T^{d_v},E_{f_2}/K)}{L_v(T^{d_v},E_{f_1}/K)} \bmod 2.
	$$
	We show that if $v \in U_{f_1,f_2}$, then $L_v(T,E_{f_1}/K) \equiv L_v(T,E_{f_2}/K) \bmod 2$. If $v \in A_{f_1} \cap A_{f_2}$, then the Euler factor is $L_v(T,E_{f_i}/K) = 1$ for $i=\{1,2\}$. Now, we can assume that $f_1 =f$ and $f_2 = 1$.
	
	Note that $M \cap M_f = M_{\text{unr}}$ and $\mid U\mid  \cap \mid U_f\mid  = \mid U\mid _{\text{unr}}$ by Lemma \ref{LemmaKodairaSymbols} (a1). If $v \in M \cap M_f$, then by Lemma \ref{LemmaLegendreSymbol} (a2), we have $L_v(T,E_f/K) \equiv L_v(T,E/K) \bmod 2$.
	Finally, if $v \in \mid U\mid  \cap \mid U_f\mid $, then by Lemma \ref{LemmaLegendreSymbol}, we have $a_v + a_{f,v} \equiv 0 \bmod 2$ and then from Equation (\ref{EulerFactorLFunction}) it follows that $L_v(T,E_f/K) \equiv L_v(T,E/K) \bmod 2$.\qedhere
\end{proof}
\section{Function Fields of Genus $0$}\label{FunctionFieldofGenus0}
In this section, we illustrate some simplifications that occur in the statements of Theorems \ref{Thm1} and \ref{Thm2} when the curve $C/k$ has genus $0$.
For simplicity, we restrict ourselves to the following context.
\begin{context} \label{context}
	Let $K := k(t)$ be a rational function field and let $E/K$ be an elliptic curve with nonconstant $j$-invariant, semistable reduction and in particular either good or multiplicative reduction over $\infty$. Suppose that over the affine plane $\mathbb{P}^1_k \smallsetminus \{\infty\}$, $E/K$ is given by a minimal Weierstrass equation of the form $y^2 = u(x)$, with $u(x) := x^3 + ax + b$ irreducible in $K[x]$ and $a,b \in k[t]$. Let $\Delta$ be the discriminant of this Weierstrass equation, which is also the discriminant of the cubic polynomial $u(x)$. Finally, we let $f \in k[t]$ be a monic square-free polynomial  of degree $d_f$ and let $S$ be the set of distinct monic irreducible polynomials in $k[t]$ dividing $f$.
\end{context}
To ease the exposition, if $K'$ is the function field of a curve $C'/k'$ over a finite field $k'$, we write $Z(T,K')$ for the zeta function $Z(T,C'/k')$. Let us first examine Theorem \ref{Thm1} in this context. For $i \in \{0,1\}$, let $L(q^iT,K_f)$ be the numerator of the zeta function $Z(q^iT, K_f)$ (see Definition \ref{DefinitionZetaFunction}).
\begin{proposition}\label{SimplificationsThm1}
	The following holds.
	\begin{enumerate}
		\item[(i)] We have the identity $L(T,\chi)L(qT,\chi) = L(T,K_f)L(qT,K_f)$.
		\item[(ii)] The set $\mid C\mid _{\text{split}} - \{\infty\}$ is
		$$
		\{ v \in k[t] \smallsetminus S \text{ monic irreducible}: x^2 - f \in k(t)[x] \text{ has a root in } k[t]/(v) \}
		$$
		and the set $\mid C\mid _{\text{inert}} - \{\infty\}$ is
		$$
		\{ v \in k[t] \smallsetminus S \text{ monic irreducible}: x^2 - f \in k(t)[x] \text{ has no roots in } k[t]/(v) \}.
		$$
		\item[(iii)] The set $\mid C\mid _{\text{ram}} - \{\infty\}$ is $S$.
		\item[(iv)] If $d_f$ is even, then $\infty \in \mid C\mid _{\text{unr}}$. If $d_f$ is odd, then $\infty \in \mid C\mid _{\text{ram}} \cap  A_{f,\text{ram}}$.
	\end{enumerate}
\end{proposition}
\begin{proof}
	\begin{enumerate}
		\item[(i)] Since the curve $C/k$ has genus $0$, the numerator of its zeta function equals $1$. Moreover, for $i \in \{0,1\}$ the denominators of $Z(q^iT,K)$ and $Z(q^iT,K_f)$ are both equal to $(1-q^iT)(1-q^{1+i}T)$. Therefore, the identity $Z(q^iT,K_f) = Z(q^iT,K)L(q^iT,\chi)$ (see \cite[p.524, (10.5) Corollary]{neukirch_1999}) simplifies to $L(q^iT,\chi) = L(q^iT,K_f)$.
		\item[(ii)] A finite place of $K$ corresponds to a monic irreducible polynomial in $k[t]$. Let $v$ be a finite place $K$ which is unramified in $K_f$. By the Dedekind-Kummer theorem (see \cite[I.3.8, pp.47-48]{neukirch_1999}), $v$ splits (resp. is inert) in $K_f$ if and only if the polynomial $x^2 - f$ has a root (resp. no roots) in the residue field $k_v = k[t]/(v)$.
		\item[(iii)] A finite place $v$ of $K$ which is ramified in $K_f$ is characterized by the fact that $\ord_v(f)$ is odd. The valuation $\ord_v(f)$ is positive (and in fact equals $1$) precisely for the finite places $v$ which correspond to monic irreducible polynomials dividing $f$.
		\item[(iv)] The place at infinity $\infty$ of $k(t)$ corresponds to the maximal ideal $(s)$ of $k[s]$ if we set $s := 1/t$. Let $w_\infty$ be a place of $K_f$ lying over the place $\infty$. Then
		$$
		\ord_{w_\infty}\left( f(t) \right) = \frac{2}{\gcd(d_f,2)}\ord_{1/t}\left( f(t) \right)
		$$
		with $2/\gcd(d_f,2)$ being the ramification index $e(w_\infty \mid \infty)$ (see \cite[p.125, Example 3.7.6]{stichtenoth_2009}). So if $d_f$ is even, then $\infty$ is unramified in $K_f$, while if $d_f$ is odd, then $t=\infty$ ramifies in $K_f$ and by Lemma \ref{LemmaKodairaSymbols} (b) it becomes a place of additive reduction for $E_f/K$.
	\end{enumerate}
\end{proof}
We now consider Theorem \ref{Thm2} in Context \ref{context}.
\begin{corollary}[of Theorem \ref{Thm2}] \label{Cor2}
	We have
	\begin{align*}
		\frac{L(T,E_f/K)}{L(T,E/K)} &\equiv(1-a_\infty T + T^2)^{\varepsilon_U}(1-T)^{\varepsilon_M}\\
		&\times \prod_{v \mid f}\left(\prod_{v \in \mid U\mid } (1 - a_vT^{d_v} + T^{2d_v}) \times \prod_{v \in M} (1-T^{d_v})\right) \bmod 2,
	\end{align*}
	
	$$
	\text{where } (\varepsilon_U,\varepsilon_M)=
	\begin{cases}
		(0,0) & \text{ if } \deg(f) \text{ is even },\\
		(1,0) & \text{ if } \deg(f) \text{ is odd and } \infty \in \mid U\mid ,\\
		(0,1) & \text{ if } \deg(f) \text{ is odd and } \infty \in M.\\
	\end{cases}
	$$
	In particular, the degree of $L(T,E_f/K)$ is equal to
	$$
	\deg L(T,E_f/K) = 2\varepsilon_U + \varepsilon_M + \sum_{v \mid f}\left( 2 \sum_{v \in \mid U \mid} d_v + \sum_{v \in M}d_v \right) + \deg(M) + 2\deg(A) - 4.
	$$
\end{corollary}
\begin{proof}
	This follows directly from Theorem \ref{Thm2}  with $f_1 = f$, $f_2 = 1$ and the following observations. If $v \in U \cup M$ is a finite place which divides $f$, then $v \in A_f$ by Lemma \ref{LemmaKodairaSymbols} (b). If $d_f$ is odd, then $\infty \in A_f$ as we just observed. Otherwise, a place of $K$ which is unramified in $K_f$ has the same reduction type before and after twisting by $f$ by Lemma \ref{LemmaKodairaSymbols} (a). The formula for the degree of $L(T,E_f/K)$ follows immediately from the expression of $L(T,E_f/K) \bmod 2$ and Equation (\ref{DegreeLFunction}).
\end{proof}
\section{Applications}\label{Application}
In this section, we discuss different applications of our results of Section \ref{ArbitraryGenus} and Section \ref{FunctionFieldofGenus0}. In particular, we prove our third main theorem (see Theorem \ref{Thm3}).

\subsection{Remark on Computational Effort} \label{Complexity}

Let $E/k(t)$ be an elliptic curve with nonconstant $j$-invariant and $L(T,E/k(t)) = 1$. If $f \in k[t]$ is a monic irreducible polynomial of even degree $d_f$ and coprime with $\Delta$, then $\deg L(T,E_f/K) = 2 d_f$. Our goal is to compare the relative costs of computing $L(T,E_f/K) \bmod 2$ using Corollary \ref{Cor2} and using the algorithm presented in \cite{baig_hall_2012}. To avoid any ambiguity, let $w$ denote the place corresponding to $f$. In this case, Corollary \ref{Cor2} gives
$$
L(T,E_f/K) \equiv 1 - a_w T^{d_f} + T^{2d_f} \bmod 2.
$$

We need to compute the single term $a_w \bmod 2$. To do this, we can proceed as follows. Since $f$ is coprime with $\Delta$, then $w \in \mid U\mid $. The Weierstrass equation $y^2 = x^3 -27 c_4 x -54 c_6$ that we gave above is minimal for the place $w$. Let $y^2 = x^3 - (27 c_4 \bmod \pi_w) x - (54 c_6 \bmod \pi_w)$ be the corresponding Weierstrass equation for the reduced elliptic curve $E_w/k_w$. Since $q$ is odd, then
$$
a_w = 1 + q^{d_f} - \#E_w(k_w) \equiv \#E_w(k_w) \bmod 2.
$$
A non-identity point of order $2$ in $E_w(k_w)[2]$, if it exists, is of the form $(\alpha,0)$, where $\alpha \in k_w$ is a root of $x^3 - (27 c_4 \bmod \pi_w)x - (54 c_6 \bmod \pi_w)$. Therefore, $x^3 - (27 c_4 \bmod \pi_w)x - (54 c_6 \bmod \pi_w)$ is irreducible in $k_w[x]$ if and only if $a_w \equiv 1 \bmod 2$.

Alternatively, we can adapt the algorithm of Baig and Hall (see \cite[pp.364-365]{baig_hall_2012}) to compute $L(T,E_f/K) \bmod 2$ as follows. One can write $L(T,E_f/K) = \sum_{n=0}^{2d_f}c_n T^n$ with $c_0 = 1$ and, for $1 \leq n \leq 2d_f$, the coefficient $c_n$ satisfies the recurrence relation (see \cite[pp.364-365, Lemma 2.2]{baig_hall_2012})
$$
c_n = \frac{1}{n} \sum_{m=1}^n b_m \cdot c_{n-m},
$$
where
$$
b_m = \sum_{\substack{v \in \mid U_f \mid \\ d_v \mid m}} d_v \cdot a_{v^{m/d_v}} + \sum_{\substack{v \in M_f^{\text{sp}} \\ d_v \mid m}} d_v + \sum_{\substack{v \in M_f^{\text{ns}} \\ d_v \mid m}} (-1)^{md_v}d_v,
$$
with the meaning, as mentioned earlier, that the subscript $f$ refers to the sets of places of $K$ over which $E_f$ has good, resp. split multiplicative, resp. nonsplit multiplicative reduction.
The amount of effort to compute $L(T,E_f/K) \bmod 2$ using this algorithm is bounded below by the computation, for each $1 \leq n \leq 2d_f$, of $a_{v^m/d_v} \bmod 2$, for $1 \leq m \leq n$ and $d_v \mid m$, where $v \in \mid U_f \mid$.
One can use the functional equation to roughly cut the number of computations in half (see \cite[p.365]{baig_hall_2012}):
$$
c_n = \varepsilon(E_f/K) q^{2n-2d_f} c_{2d_f-n},~~0 \leq n \leq 2 d_f.
$$
For each degree $d_v$, there are $\frac{1}{d_v}\sum_{ d \mid d_v } \mu(d_v \mid d) \cdot q^{d_v}$ monic irreducible polynomials of degree $d_v$ over the finite field $k_v$ with $q^{d_v}$ elements (see \cite[pp.567-568, Chapter 14]{dummit_foote_2004}), where $\mu$ is the M\"obius function. Most of these polynomials correspond to places of good reduction for $E_f/K$ and we roughly need to compute half of the quantities of the form $a_{v^m/d_v} \bmod 2$ that appears in the formulas for the $b_m$'s. Our formula is therefore clearly more efficient to compute $L(T,E_f/K) \bmod 2$ than the algorithm in \cite{baig_hall_2012}.
\subsection{Computing $L$-Functions of Elliptic Curves \`a la Schoof}\label{ComputationLFunction}

In this subsection, we discuss a general strategy to compute polynomial $L$-functions of elliptic curves. Our strategy is motivated by Schoof's algorithm which computes zeta functions of elliptic curves over finite fields. We also prove the elementary Proposition \ref{OddAnalyticRank} on elliptic curves of odd analytic rank that could improve our approach. We illustrate this strategy in Proposition \ref{HallDegree2} and in Proposition \ref{ComputationQuadraticTwists} by computing some degree $2$ $L$-functions using Theorem \ref{HallThm}, Theorem \ref{Thm1}, Corollary \ref{Cor2} and Proposition \ref{OddAnalyticRank}.

Recall from Subsection \ref{RecallsLFunctions} that the $L$-function of $E/K$ is a polynomial $L(T,E/K) = \sum_{n=0}^d a_nT^n$ in $1 + T \cdot \Z[T]$ of degree $d = 2(2g(C)-2) + \deg(M) + 2\deg(A)$ and the leading coefficient $\varepsilon(E/K)q^d$. We can write $L(T,E/K) = \prod_{i=1}^d(1 - \alpha_n T)$ as product of its reciprocal roots
$\alpha_1, \ldots, \alpha_d$. The reciprocal roots are algebraic integers which satisfy $\mid \alpha_i \mid = q$, for each $1 \leq i \leq d$ and each complex embedding (see \cite[p.267, Theorem 2.2.1]{ulmer_2011}). Using this fact and expressing the coefficients of $L(T,E/K)$ as symmetric polynomials in the $\alpha_i$'s, we deduce that
$$
\mid a_n\mid  \leq \binom{d}{n}q^n \text{ for each } 1 \leq n \leq d.
$$
Given an explicit Weierstrass equation for $E/K$ we can already compute $a_d = \varepsilon(E/K)q^d$ (using \cite[p.365, Lemma 2.3]{baig_hall_2012}). Therefore, we only need to compute the coefficients $a_i, 1 \leq i \leq d-1$. To guarantee in general that $L(T,E/K)$ is completely determined by its reduction $L(T,E/K) \bmod M$ for some positive integer $M$ coprime to $q$, it suffices that $M$ satisfies
$$
M > 2\max_{1 \leq n \leq d-1} \binom{d}{n}q^n.
$$
\begin{remark}\label{ObservationGlobalRootNumber}
	If $E(K)$ has a subgroup of order $N \geq 3$, with $N$ coprime to $q$, then we can avoid the use of Tate's algorithm (which is needed in \cite[p.365, Lemma 2.3]{baig_hall_2012}), and compute global root numbers more directly. Indeed, expanding the formula for $L(T,E_f/K) \bmod N$ in Theorem \ref{Thm1} gives a relation of the form $\varepsilon(E_f/K)q^d \equiv \pm q^d \bmod N$ and so $\varepsilon(E_f/K) \equiv \pm 1 \bmod N$. This is sufficient to determine $\varepsilon(E_f/K)$ since there is only one representative for the class $\varepsilon(E_f/K) \bmod N$ in the set $\{-1,1\}$. Alternatively, one could use \cite[p.134, Corollary 5]{hall_2006} and Proposition \ref{DegreeAndRootNumberQuadraticTwists}.
\end{remark}
The following simple observation does not seem to be in the literature.
\begin{proposition}\label{OddAnalyticRank}
	If $\varepsilon(E/K) = -1$ and $d = \deg(L(T,E/K))$, then
	\begin{equation}\label{EquationOddAnalyticRank}
		L(T,E/K) = (1-qT)\left( 1 + \sum_{n=1}^{d-2}\left( \sum_{m=0}^n a_m q^{n-m} \right)T^n + q^{d-1}T^{d-1} \right). 
	\end{equation}
\end{proposition}
\begin{proof}
	On one hand, since $E/K$ has odd analytic rank, then
	$$
	0 = L(q^{-1},E/K) = \sum_{n=1}^{d-1}a_n(q^{-1})^n \text{ and so } a_{d-1} = -\sum_{n=1}^{d-2}a_n q^{d-(n+1)}.
	$$
	On the other hand,	expanding the right hand side of Equation \eqref{EquationOddAnalyticRank}, gives
	\begin{align*}
		&(1-qT)\left( 1 + \sum_{n=1}^{d-2}\left( \sum_{m=0}^n a_m q^{n-m} \right)T^n + q^{d-1}T^{d-1} \right)\\
		&= 1 + \sum_{n=1}^{d-2}a_nT^{d-2} -\left(\sum_{n=1}^{d-2}a_n q^{d-(n+1)}\right)T^{d-1} - q^dT^d.
	\end{align*}
\end{proof}
In particular, we have the following consequence.
\begin{corollary}\label{FirstCasesOddAnalyticRank}
	Under the assumptions of Proposition \ref{OddAnalyticRank} we have
	\begin{enumerate}
		\item[(i)] If $d = 1$, then $L(T,E/K) = 1-qT$.
		\item[(ii)] If $d = 2$, then $L(T,E/K) = (1-qT)(1+qT)$.
		\item[(iii)] If $d = 3$, then $L(T,E/K) = (1-qT)\left( 1 + (q+a_1)T + q^2T^2\right)$ and $a_2 = -7a_1$.
	\end{enumerate}
\end{corollary}
\begin{remark}\label{ComputationalRemark}
	A priori, if $d = 3$, then we have $\mid a_1\mid  \leq 3q$ and $\mid a_2\mid  \leq 3q^2$. However, if $E/K$ has odd analytic rank, then Corollary \ref{FirstCasesOddAnalyticRank} (iii) shows that we only need to consider the bound $\mid a_1\mid  \leq 3q$ in order the compute $L(T,E/K)$. More generally, Proposition \ref{OddAnalyticRank} shows that if $E/K$ has odd analytic rank and $L(T,E/K)$ has degree $d$, then in order to compute $L(T,E/K)$, we only need to consider the bounds on $\mid a_i\mid $, for $1 \leq i \leq d - 2$.
\end{remark}
\begin{proposition}\label{HallDegree2}
	Let $k$ be the field with $5$ elements and let $f$ be a rational function of the form $f \in \{ P(t), 1/Q(t), P(t)/Q(t), \text{ with } P(t), Q(t) \in k[t] \text{ coprime and of degree } 1 \}$. An elliptic curve $E/k(t)$ given by  a Weierstrass equation
	$$
	E/k(t) : y^2 = (1-a)xy - by = x^3 - bx^2,
	$$
	with $a = (10-2f)/(f^2-9)$ and $b = -2(f-1)^2(f-5)/(f^2-9)^2$, has torsion subgroup $E(k(t))_\tors$ isomorphic to $\Z/2\Z \times \Z/6\Z$. Moreover, its $L$-function equals
	$$
	L(T,E/k(t)) = (1-5T)(1+5T).
	$$
\end{proposition}
\begin{proof}
	From \cite[p.403, Table 3]{mcdonald_2018}, the prime-to-$5$ torsion subgroup of $E(k(t))$ is isomorphic to $\Z/2\Z \times \Z/6\Z$ and since the $j$-invariant of $E$,
	$$
	j(E) = 4\frac{(f^2+f+1)^3(f^2+2f+3)^3(f^2+4f+1)^3(f^2+4f+2)^3}{f^6(f+1)^2(f+2)^2(f+3)^2(f+4)^6},
	$$
	is not a $5$th power in $k(t)$, then the $5$-torsion subgroup of $E(k(t))$ is trivial (see \cite[p.229, Proposition 7.3]{ulmer_2011}). By direct computation, we see that the places of $k(t)$ over which $E$ has bad reduction is $M$ and that $\#M^{\text{sp}} = \#M^{\text{ns}} = 3$ and so $\deg(L(T,E/k(t))) = 2$. Since $k(t)$ has genus $0$, then $Z(T,k(t))Z(5T,k(t)) = 1/( (1-T)(1-5T)^2(1-5^2T) )$. Now, $5^2 \equiv 1 \bmod 12$ and that we have the factorizations $(1-5T)^2 \equiv (1+T)^2 \bmod 12$ and $(1-5T)(1-T) \equiv (1+T)(1+5T) \bmod 12$. It follows from Theorem \ref{HallThm} that
	\begin{align*}
		L(T,E/k(t)) &\equiv \frac{1}{(1-T)^2(1-5T)^2} \times (1-5T)^3 \times \frac{(1-T)^3(1-5T)^3}{(1+T)^3}\\
		&\equiv \frac{(1-T)(1-5T)^4}{(1+T)^3}\\
		&\equiv (1-5T)(1+5T) \bmod 12.
	\end{align*}
	Then, $\varepsilon(E/k(t))5^2 \equiv -5^2 \bmod 12$ and so $\varepsilon(E/k(t)) = -1$ (see Remark \ref{ObservationGlobalRootNumber}). To conclude that $L(T,E/k(t)) = (1-5T)(1+5T)$, we can either use Corollary \ref{FirstCasesOddAnalyticRank} (ii), or observe that since $a_1 \equiv 0 \bmod 12$ and $\mid a_1\mid  \leq 10$, then $a_1 = 0$.
\end{proof}

\begin{proposition}\label{ComputationQuadraticTwists}
	Suppose that $k$ is the field with $7$ elements and let $f \in k[t]$ be a polynomial of degree $1$. An elliptic curve $E/k(t)$ given by a Weierstrass equation
	$$
	E/k(t): y^2 + xy + \frac{f^2+f+1}{3(f+2)^3} = x^3
	$$
	has torsion subgroup $E(k(t))_\tors$ isomorphic to $(\Z/3\Z)^2$, bad reduction $M^{\text{sp}}$ and $L$-function equal to $1$. Moreover, for any finite place $h \in M^{\text{sp}}$, the $L$-function of the quadratic twist $E_h/k(t)$ is
	$$
	L(T,E_h/k(t)) = (1-7T)(1+7T).
	$$
\end{proposition}
\begin{proof}
	Since $char(k) \neq 3$ and $\zeta_3 \in k$, then the prime-to-$7$ torsion subgroup of $E(k(t))$ is isomorphic to $(\Z/3\Z)^2$ (see \cite[p.402, Table 1]{mcdonald_2018}). Since the $j$-invariant of $E$,
	$$
	j(E) = 6\frac{f^3(f+1)^2(f+2)^3(f+4)^2}{(f+6)^3},
	$$
	is not a $7$th power in $k(t)$, then the $7$-torsion subgroup of $E(k(t))$ is trivial (see \cite[p.229, Proposition 7.3]{ulmer_2011}). By direct computation, we see that the set of places of $k(t)$ over which $E$ has bad reduction is $M^\text{sp}$ and the elements of this set are $3$ finite places of degree $1$ (none of which is a monic irreducible factor of $f$) and the place $\infty$. In particular, $\deg(L(T,E/k(t))) = 0$.
	
	We now consider the quadratic twists $E_h/k(t)$ for any finite place $h \in M^{\text{sp}}$. Using Lemma \ref{SimplificationsThm1}, we see that $\# M_\text{split}^{sp} = \# M_\text{inert}^{sp} = 1$ and $M_\text{ram}^{sp} = \{h,\infty\}$. From Lemma \ref{LemmaKodairaSymbols} (a2), we have $M_\text{split}^{sp}= M_h^{\text{sp}}$ and $M_\text{inert}^{sp} = M_h^{\text{ns}}$ and from Lemma \ref{SimplificationsThm1} (iv), we have $A_h = \{h, \infty \}$. So $\deg(L(T,E_h/k(t)) = -4 + 2 + 2(2) = 2$. By the Riemann-Hurwitz formula, we see that the quadratic extension $k(t)_h/k(t)$ has genus $0$ and so $L(T,\chi)L(7T,\chi) = 1$.  Therefore, Theorem \ref{Thm1} implies that
	$$
	L(T,E_h/k(t)) \equiv (1-7T)(1+7T) \bmod 9.
	$$
	Then, $\varepsilon(E_h/k(t))7^2 \equiv -7^2 \bmod 9$ and so $\varepsilon(E_h/k(t)) = -1$ (see Remark \ref{ObservationGlobalRootNumber}). To conclude that $L(T,E_h/k(t)) = (1-7T)(1+7T)$, we can either use Corollary \ref{FirstCasesOddAnalyticRank} (ii) or proceed as follows. We have $a_1 \equiv 0 \bmod 9$ and since $L(T,E/k(t)) = 1$, $\deg(h) = 1$ and $\infty \in M$, then Corollary \ref{Cor2} implies that
	$$
	L(T,E_h/k(t)) \equiv (1-T)^2 \bmod 2.
	$$
	Hence, $a_1 \equiv 0 \bmod 2$ and by the Chinese remainder theorem, $a_1 \equiv 0 \bmod 18$. Since $\mid a_1\mid  \leq 14$, we have $a_1 = 0$.
\end{proof}
\subsection{Study of an Infinite Family of Quadratic Twists}\label{StudyInfiniteFamily}
In this subsection, we apply Theorem \ref{Thm1}, Corollary \ref{Cor2}, Lemma \ref{NumeratorZetaFunctionAt-1} and Proposition \ref{DegreeAndRootNumberQuadraticTwists} to study an infinite family of quadratic twists of the so-called ``universal'' elliptic curve $E$ over the function field $K=k(t)$ of the ``modular curve $X_1(3)/k$''. We suppose that $p \geq 5$. The elliptic curve $E/K$ is defined by the Weierstrass equation
\begin{equation}\label{ThreeTorsion}
	y^2 + 3xy + (1-t^3)y = x^3
\end{equation}
together with the point at infinity $O = [0:1:0]$. Its discriminant is $\Delta = -3^3t^3(t-1)^3(t^2+t+1)^3$.
Since $char(k) \neq 2,3$, we can rewrite Equation (\ref{ThreeTorsion}) into the Weierstrass equation
$$
y^2 = x^3 - 27c_4 x - 54c_6, 
$$
with $c_4 = 3^2(2t+1)(4t^2-2t+1)$ and $c_6 =-3^3(2t^2+2t-1)(4t^4-4t^3+6t^2+2t+1)$ (see \cite[p.43]{silverman_2009}).
The $j$-invariant of $E/K$ is the nonconstant rational function
$$
j(E) = -\frac{3^3(2t+1)^3(4t^2-2t+1)^3}{t^3(t-1)^3(t^2+t+1)^3}
$$
and the point $P=(0,0)$ has order $3$ since $[-1]P = (0,1-t^3) = [2]P$ (see \cite[p.54, (a),(c)]{silverman_2009}).

The discriminant $\Delta$ is minimal for all the finite places of $E/K$. In particular, the finite places of bad reduction of $E/K$ correspond to $t,t-1$ and the irreducible factors of $t^2+t+1$. The latter splits into two linear factors if and only if its discriminant $-3$ is a square in $k$. In terms of Legendre symbol, we have
\begin{equation}\label{LegendreSymbol}
	\left( \frac{-3}{q} \right) =
	\begin{cases}
		1 & \text{ if } q \equiv 1, 7 \bmod 12\\
		-1 & \text{ if } q \equiv 5, 11 \bmod 12.\\
	\end{cases}
\end{equation}
The elliptic curve $E/K$ has good reduction above the place $\infty$. Indeed, there is a smallest integer $e \geq 1$ for which applying the change of variables $(t,x,y) \mapsto (1/s,x/s^{2e}, y/s^{3e})$ to Equation (\ref{ThreeTorsion}) yields a Weierstrass equation over $k[s]$ that is minimal at $s$ : $e=1$. The resulting Weierstrass equation is
$$
y^2 = x^3 -27 c_{4,\infty} x - 54c_{4,\infty}, \text{where } c_{4,\infty} = 3^2(s^3+8) \text{ and } c_{6,\infty} = 3^3(s^6-20s^3-8).
$$
Since $\ord_\infty(c_{4,\infty}) = 0$, this Weierstrass equation is minimal over $\infty$ and its discriminant $\Delta_\infty = 3^2(s-1)^3(s^2+s+1)^3$ satisfies $\ord_\infty(\Delta_\infty) = 0$. Therefore, $E/K$ has good reduction over $\infty$. In particular, $\deg(L(T,E/K)) = 0$.

We now assume that $q \equiv 5,11 \bmod 12$, so that the polynomial $t^2+t+1$ is irreducible in $k[t]$. We prove the following result.

\begin{theorem}\label{Thm3}
	Let $E/K$ be the universal elliptic curve over the function field $K = k(t)$ of the modular curve $X_1(3)/k$, where $k$ is a finite field of cardinality $q$. Suppose that $q \equiv 5, 11 \bmod 12$. Let $f \in k[t]$ be a monic square-free polynomial of degree $d_f$ which is coprime to $\Delta$. Then
	$$
	\deg\left( L(T,E_f/K) \right) =
	\begin{cases}
		2d_f &\text{ if } d_f \text{ is even},\\
		2(d_f + 1) &\text{ if } d_f \text{ is odd}
	\end{cases}
	$$
	is an even number. For $a \in \{0,1\}$, we write $f(a) \in k^2$? to ask whether $f(a)$ is a square (Yes) or not (No) in $k$. We have

		\begin{table}[h]
			\centering
		\begin{tabular}{|c|c|c|}
			\hline
			$f(0) \in k^2$? & $f(1) \in k^2$? & $L(T,E_f/K) \bmod 3$\\
			\hline
			Yes & Yes & $L(T,K_f)L(-T,K_f)(1-T)(1+T)(1+T^2)$\\
			\hline
			Yes & No & $L(T,K_f)L(-T,K_f)(1-T)^2(1+T^2)$\\
			\hline
			No & Yes & $L(T,K_f)L(-T,K_f)(1-T)(1+T)^3$\\
			\hline
			No & No & $L(T,K_f)L(-T,K_f)(1-T)(1+T)(1+T^2)$\\
			\hline
		\end{tabular}.
		\caption{Reduction $L(T,E_f/K) \bmod 3$}
	\label{tb14ReductionMod3}
	\end{table}

	The global root number $\varepsilon(E_f/K)$ is given in Table \ref{tb15GlobalRootNumberandAnalyticRank} below. 

	\begin{table}[h]
		\centering
		\begin{tabular}{ |c|c|c|c|}
			\hline
			$f(0) \in k^2$? & $f(1) \in k^2$? &  $\varepsilon(E_f/K)$ & $\mathrm{ran}^{an}(E_f/K)$\\
			\hline
			Yes & Yes & $1$ & $0$\\
			\hline
			Yes & No & $1$ & $0$ \\
			\hline
			No & Yes & $-1$ & $\leq 3$ \\
			\hline
			No & No & $-1$ &  $1$\\
			\hline
		\end{tabular}
		\caption{Global Root Number and Analytic Rank Using Reduction Modulo 3}
		\label{tb15GlobalRootNumberandAnalyticRank}
	\end{table}

	Moreover, let $\mathrm{Jac}(C_f)(k)$ be the set of $k$-rational points of the Jacobian of the (normalization of the) curve $C_f/k$ associated to the field $K_f$. If the $3$-torsion subgroup of $\mathrm{Jac}(C_f)(k)$ is trivial, then the fourth column of Table \ref{tb15GlobalRootNumberandAnalyticRank} gives the analytic rank, $\mathrm{rank}^{an}(E_f/K)$, of the elliptic curve $E_f/K$ or an upper bound on it.
	
	We also have
	$$
	L(T,E_f/K) \equiv (1-a_\infty T + T^2)^{\varepsilon_U} \times \prod_{v \mid f}(1 - a_vT^{d_v} + T^{2d_v}) \bmod 2,
	$$
	where $\varepsilon_U = 1$ if $d_f$ is odd and $\varepsilon_U = 0$ if $d_f$ is even. If $d_f$ is odd, then
	$$
	\mathrm{rank}^{an}(E_f/K) \leq
	\begin{cases}
		2\# \{ v \mid f : a_v \equiv 0 \bmod 2 \} + 2 &\text{ if } a_\infty \equiv 0 \bmod 2,\\
		2\# \{ v \mid f : a_v \equiv 0 \bmod 2 \} &\text{ if } a_\infty \equiv 1 \bmod 2.\\
	\end{cases}
	$$
	In particular, if $f$ is irreducible of odd degree, then

	\begin{table}[h]
		\centering
		\begin{tabular}{ |c|c|c|}
			\hline
			$a_\infty \bmod 2$ & $a_f \bmod 2$ & $\mathrm{rank}^{an}(E_f/K)$ \\
			\hline
			$0$ & $0$ & $\leq 4$\\
			\hline
			$0$ & $1$ & $\leq 2$\\
			\hline
			$1$ & $0$ & $\leq 2$\\
			\hline
			$1$ & $1$ & $ 0$\\
			\hline
		\end{tabular}
		\caption{Analytic Rank using Reduction Modulo $2$}
	\label{tb16Mod2}
	\end{table}

	If $f = \Delta$, then $\deg(L(T,E_\Delta/K)) = 4$, $\varepsilon(E_\Delta/K) = 1$ and
	$$
	L(T,E_\Delta/K) \equiv
	\begin{cases}
		L(T,K_f)L(-T,K_f) \bmod 3,\\
		(1-T)^4 \bmod 2.
	\end{cases}
	$$
\end{theorem}
\begin{proof}
	We first qualify the places of $K$ over which $E$ has bad reduction. The Weierstrass equation $y^2 = x^3 -27 c_4 x - 54c_6$ is minimal for $v \in S := \{ (t), (t-1), (t^2+t+1) \}$ since $-27c_4, -54c_6 \in k[t]$ and for each $v \in S, \ord_v(\Delta) = 3 < 12$ (see \cite[Remark 1.1, p.186]{silverman_2009}). Moreover, for $v \in S$, $\ord_v(c_4) = 0$ and so $M = S$ with $\deg(M) = 4$ and $A = \emptyset$.
	From Equation (\ref{DegreeLFunction}), $\deg(L(T,E/K)) = 0$ and thus $L(T,E/K) = 1$. By \cite[p.364, Lemma 2.1]{baig_hall_2012}, a place $v \in M$ is in $M^{\text{sp}}$ if and only if the image of $6(-54c_6)$ in the residue field $k_v$ is a square. We have $6(-54c_6) \equiv -2^2 \cdot 3^7 \mod \pi_t$. This is not a square in $k_t$ since $-3$ is not a square in $k_t$ as we assumed that $q \equiv 5, 11 \bmod 12$ (see Equation (\ref{LegendreSymbol})). Now, if $v \in \{ (t-1), (t^2+t+1) \}$, then $6(-54c_6) \equiv (2^2 \cdot 3^5)^2 \bmod \pi_v$. Therefore, $M^{\text{sp}} = \{ (t-1), (t^2+t+1) \}$ and $M^{\text{ns}} = \{ (t)\}$.
	
	We now study the quadratic twists $E_f/K$ and their $L$-functions and first consider the reduction modulo $3$. Suppose until further notice that $f \in k[t]$ is a monic square-free polynomial of degree $d_f$ that is coprime to $\Delta$. By Proposition \ref{DegreeAndRootNumberQuadraticTwists} (i), we have
	$$
	\deg\left( L(T,E_f/K) \right) =
	\begin{cases}
		2d_f &\text{ if } d_f \text{ is even},\\
		2(d_f + 1) &\text{ if } d_f \text{ is odd},
	\end{cases}
	$$
	which is in particular an even number. Moreover $M^\text{sp} = M_\unr^{\text{sp}}$ and it follows from Lemma \ref{SimplificationsThm1} that $M_\inert^{\text{ns}} = \{t\}$ if and only if $f(0)$ is not a square in $k$. Therefore, from Proposition \ref{DegreeAndRootNumberQuadraticTwists} (ii) we have
	$$
	\varepsilon(E_f/K) =
	\begin{cases}
		-1 & \text{ if } f(0) \in k \smallsetminus k^2,\\
		1 & \text{ if } f(0) \in  k^2.
	\end{cases}
	$$
	Since $q \equiv 5,11 \bmod 12$, then $q \equiv - 1 \bmod 3$. In this situation, Theorem \ref{Thm1} gives
	$$
	L(T,E_f/K) \equiv L(T,K_f)L(-T,K_f)(1+\alpha_tT)(1 - \alpha_{t-1}T)(1+\alpha_{t^2+t+1}T^2) \bmod 3,
	$$
	where for $v \in M$, $\alpha_v = 1$ if $v \in M_{\inert}$ and $\alpha_v = -1$ if $v \in M_{\text{split}}$. From Lemma \ref{SimplificationsThm1}, we see that
	\begin{equation}\label{InertSplitX13}
		\alpha_t =
		\begin{cases}
			-1 & \text{ if } f(0) \in k^2,\\
			1 & \text{ if } f(0) \in k \smallsetminus k^2
		\end{cases}
		\text{ and }
		\alpha_{t-1} =
		\begin{cases}
			-1 & \text{ if } f(1) \in k^2,\\
			1 & \text{ if } f(1) \in k \smallsetminus k^2.
		\end{cases}
	\end{equation}
	Now, we have $L(T,K_f)L(-T,K_f) = q^{2g(K_f)}T^{4g(K_f)} + O( T^{4g(K_f)-1})$ (see \cite[p.54, Proof of Theorem 5.9]{rosen_2002}) and therefore we have
	$$
	L(T,E_f/K) \equiv -\alpha_t \alpha_{t-1} \alpha_{t^2+t+1} T^{4(g(K_f)+1)} + O(T^{4g(K_f) + 3}) \bmod 3
	$$
	and in particular
	\begin{equation}\label{RelationGlobalRootNumberInertiaSplitX13}
		\varepsilon(E_f/K) \equiv -\alpha_t \alpha_{t-1} \alpha_{t^2+t+1} \bmod 3.
	\end{equation}
	Combining Equation \eqref{InertSplitX13} and Equation \eqref{RelationGlobalRootNumberInertiaSplitX13} we deduce the following.
	\begin{table}[h]
		\centering
		\begin{tabular}{ |c|c|c|c|}
			\hline
			$f(0) \in k^2$? & $f(1) \in k^2$? & $(\alpha_t, \alpha_{t-1}, \alpha_{t^2+t+1})$ & $\varepsilon(E_f/K)$ \\
			\hline
			Yes & Yes & $(-1,-1,-1)$ & $1$ \\
			\hline
			Yes & No & $(-1,1,1)$ & $1$ \\
			\hline
			No & Yes & $(1,-1,-1)$ & $-1$ \\
			\hline
			No & No & $(1,1,1)$ & $-1$ \\
			\hline
		\end{tabular}
		\caption{Decomposition of Unramified Places in $K_f$ and Global Root Number}
	\label{tb17UnramifiedRootNumber}
	\end{table}

	Putting all these together, we conclude that
	\begin{table}[h]
		\centering
		\begin{tabular}{|c|c|c|c|}
			\hline
			$f(0) \in k^2$? & $f(1) \in k^2$? & $L(T,E_f/K) \bmod 3$ & $\varepsilon(E_f/K)$ \\
			\hline
			Yes & Yes & $L(T,K_f)L(-T,K_f)(1-T)(1+T)(1+T^2)$ & $1$ \\
			\hline
			Yes & No & $L(T,K_f)L(-T,K_f)(1-T)^2(1+T^2)$ & $1$ \\
			\hline
			No & Yes & $L(T,K_f)L(-T,K_f)(1-T)(1+T)^3$ & $-1$ \\
			\hline
			No & No & $L(T,K_f)L(-T,K_f)(1-T)(1+T)(1+T^2)$ & $-1$ \\
			\hline
		\end{tabular}
		\caption{Reduction $L(T,E_f/K) \bmod 3$ and Global Root Number}
	\label{tb18Mod3andRootNumber}
	\end{table}
		
	Similarly to \cite[p.134, Corollary 6]{hall_2006}, if $\#E(K)$ is divisible by $3$, then we define, for any polynomial $P(T) \in 1 + T \cdot \Z[T]$, $v_3(P)$ as the $(1+T)$-adic valuation of the image of $P$ in $1 + T \cdot \F_3[T]$ and extend it linearly to rational functions $P_1(T)/P_2(T)$. We have $L(1,C_f/k) = \#\mathrm{Jac}(C_f)(k)$ (see \cite[p.53, Theorem 5.9]{rosen_2002}) and since $q \equiv -1 \bmod 3$, then $L(-1,C_f/k) \equiv (-1)^{g(C_f)} \#\mathrm{Jac}(C_f)(k) \bmod 3$, by Lemma \ref{NumeratorZetaFunctionAt-1}. Therefore, we have $\mathrm{Jac}(K_f)(k)[3] = \{O\}$ if and only if $v_3(L(T,K_f)) = v_3(L(-T,K_f)) = 0$. Hence, in general we have
	$$
	\mathrm{rank}^{an}(E_f(K)) \leq v_3(L(T,K_f)) + v_3(L(-T,K_f)) + v_3(R(T)),
	$$
	where

	\begin{table}[h]
		\centering
		\begin{tabular}{ |c|c|c|c|}
			\hline
			$f(0) \in k^2$? & $f(1) \in k^2$? & $v_3(R(T))$ & $\varepsilon(E_f/K)$ \\
			\hline
			Yes & Yes & $1$ & $1$ \\
			\hline
			Yes & No & $0$ & $1$ \\
			\hline
			No & Yes & $3$ & $-1$ \\
			\hline
			No & No & $1$ & $-1$ \\
			\hline
		\end{tabular}
		\caption{$(1+T)$-adic valuation in $\F_3[T]$ and Global Root Number}
	\label{tb19AdicValuationAndRootNumber}
	\end{table}

	and in particular, if $\mathrm{Jac}(K_f)[3] = \{O\}$, then we obtain Table \ref{tb20RootNumberAnalyticRank}.

	\begin{table}[h]
		\centering
		\begin{tabular}{ |c|c|c|c|}
			\hline
			$f(0) \in k^2$? & $f(1) \in k^2$? & $\varepsilon(E_f/K)$ & $\mathrm{ran}^{an}(E_f/K)$\\
			\hline
			Yes & Yes & $1$ & $0$ \\
			\hline
			Yes & No & $1$ & $0$ \\
			\hline
			No & Yes & $-1$ & $\leq 3$ \\
			\hline
			No & No & $-1$ & $1$ \\
			\hline
		\end{tabular}
		\caption{Global Root Number and Analytic Rank Using Reduction Modulo $3$}
	\label{tb20RootNumberAnalyticRank}
	\end{table}

	Now suppose that $f = \Delta$. The set of places of $K$ over which $E$ has bad reduction is in $M$ and then by Lemma \ref{SimplificationsThm1} (iii) these places ramify in $K_\Delta$ and acquire additive reduction with respect to $E_\Delta/K$. Therefore, $\deg\left( L(T,E_\Delta/K) \right) = -4 + 2(4) = 4$.
	By Theorem \ref{Thm1} we have
	$$
	L(T,E_\Delta/K) \equiv L(T,C_f/k)L(-T,C_f/k)\bmod 3,
	$$
	with $L(T,C_f/k) = qT^2 + O(T)$. We have the relation $\varepsilon(E_\Delta/K)q^4 \equiv q^2 \bmod 3$ and so $\varepsilon(E_\Delta/K) = 1$ (see Remark \ref{ObservationGlobalRootNumber}). If $\mathrm{Jac}(C_f)(k)[3] = \{O\}$, then $\mathrm{rank}^{an}(E_\Delta/K) = 0$.
	
	We now consider the reduction modulo $2$. The places of $K$ over which $E$ has bad reduction are all in $M$, and $\infty \in \mid U\mid $. Moreover, we know that $L(T,E/K) = 1$. Corollary \ref{Cor2} allows us to describe $L(T,E_f,K) \bmod 2$. More precisely, if $f$ is coprime with $\Delta$, then
	\begin{align*}
		L(T,E_f/K) &\equiv (1-a_\infty T + T^2)^{\varepsilon_U} \times \prod_{v \mid f}(1 - a_vT^{d_v} + T^{2d_v}) \bmod 2\\
		&\equiv (1-a_\infty T + T^2)^{\varepsilon_U} \times \prod_{v \mid f}( \prod_{a_v \equiv 0 \bmod 2} (1 - T)^2 (1 + T + \ldots + T^{d_v-1})^2\\
		&\times \prod_{a_v \equiv 1 \bmod 2} (1 + T^{d_v} + T^{2d_v})) \bmod 2.
	\end{align*}
	where
	$$
	\varepsilon_U =
	\begin{cases}
		0 & \text{ if } d_f \text{ is even,}\\
		1 & \text{ if } d_f \text{ is odd.}
	\end{cases}
	$$
	From this, we see for example that if $d_f$ is odd, then
	$$
	\mathrm{rank}(E_f/K) \leq
	\begin{cases}
		2\# \{ v \mid f : a_v \equiv 0 \bmod 2 \} + 2 &\text{ if } a_\infty \equiv 0 \bmod 2,\\
		2\# \{ v \mid f : a_v \equiv 0 \bmod 2 \} &\text{ if } a_\infty \equiv 1 \bmod 2.\\
	\end{cases}
	$$
	In particular, if $f$ is irreducible of odd degree, we have the following.

	\begin{table}[h]
		\centering
		\begin{tabular}{ |c|c|c|}
			\hline
			$a_\infty \bmod 2$ & $a_f \bmod 2$ & $\mathrm{ran}_{an}(E_f/K)$ \\
			\hline
			$0$ & $0$ & $\leq 4$\\
			\hline
			$0$ & $1$ & $\leq 2$\\
			\hline
			$1$ & $0$ & $\leq 2$\\
			\hline
			$1$ & $1$ & $= 0$\\
			\hline
		\end{tabular}
	\caption{Analytic Rank Using Reduction Modulo $2$}
	\label{tb21Mod2AnalyticRank}
\end{table}

	Finally, $f = \Delta$ has even degree and so
	\begin{equation}\label{RamifiedTwistMod2}
	L(T,E_\Delta/K) \equiv (1-T)^2(1-T^2) \equiv (1-T)^4 \bmod 2.
	\end{equation}
	Then the bound on the analytic rank provided by Equation \eqref{RamifiedTwistMod2} is $4 = \deg( E_\Delta/K )$.
\end{proof}
\section{Appendix - Proofs of Some Lemmas}\label{Appendix}
In this section we prove some lemmas which we believe are already known, at least in the folklore, but which we could not find in the literature.

Let $v$ be a place of $K$ and let $K_v$ be its $v$-adic completion. Let $\pi_v$ be a chosen normalized uniformizer. There is a well-defined reduction map $E(K_v) \to E_v(k_v)$ modulo $\pi_v$ (see \cite[p.187, VII.2]{silverman_2009}).
The subset $E_{v,\sm}(k_v)$ of smooth points of $E_v(k_v)$ forms a group (see \cite[p.56, III.2.5 and p.105, Exercise 3.5]{silverman_2009}). Now, let $E_0(K_v)$ be the subset of points of $E(K_v)$ whose reduction is smooth and let $E_1(K_v)$ be the kernel of the reduction map. There is a short exact sequence of abelian groups
\begin{equation}
	0 \to E_1(K_v) \to E_0(K_v) \to E_{v,\sm}(k_v) \to 0, \label{SesE1E0Esm}
\end{equation}
where the map $E_0(K_v) \to E_{v,\sm}(k_v)$ is the reduction modulo the uniformizer $\pi_v$ (see \cite[p.188, VII.2.1]{silverman_2009}). Our proof of the following lemma uses the language of N\'eron models (See for example \cite[pp.361-364, IV.9]{silverman_1994}.).
\begin{lemma}\label{AppendixPrimeTopTorsion}
	Let $K_v$ be a field complete with respect to a non-Archimedean place $v$ whose residue field $k_v$ has characteristic $p \geq 0$. Let $E/K_v$ be an elliptic curve with additive reduction over $v$. Then, the prime-to-$p$-torsion subgroup of $E(K_v)$ injects into the group $\Phi_v(k_v)$.
\end{lemma}
\begin{proof}
	If $N$ is an integer coprime with $p$, then $E_1(k_v)[N]=\{O\}$ (see \cite[p.192, VII.3.1(a)]{silverman_2009}) and since $E/K_v$ has additive reduction, then $E_{v,\sm}(k_v)$ is isomorphic to the additive $p$-group $k_v$ (see \cite[p.56, III.2.5 and p.105, Exercise 3.5]{silverman_2009}). Therefore, $E_{v,\sm}(k_v)[N] = \{O\}$. In view of the short exact sequence (\ref{SesE1E0Esm}), we deduce that $E_0(K_v)[N] = \{O\}$. Now, there is a short exact sequence of abelian groups
	$$
	0 \to E_0(K_v) \to E(K_v) \to E(K_v)/E_0(K_v) \to 0.
	$$
	Since $E_0(K_v)[N] = \{O\}$, then $E(k_v)[N]$ injects into $E(K_v)/E_0(K_v)$. Let $\calE_{k_v}$ (resp. $\calE_{k_v}^0$) be the special fiber of the N\'eron model $\calE_v/\Spec(\calO_v)$ of $E/K_v$ (resp. of the identity component $\calE_v^0/\Spec(\calO_v)$ of $\calE_v/\Spec(\calO_v)$). Since the field $K_v$ is complete (and in particular Henselian), then from \cite[p.362, 9.2(b)]{silverman_1994} we have the group isomorphisms
	$$
	E(K_v)/E_0(K_v) \simeq \calE_{k_v}(k_v)/\calE^0_{k_v}(k_v) \simeq \Phi_v(k_v).
	$$
	Therefore, the prime-to-$p$-part $\cup_{(N,p)=1} E(K_v)[N]$ of $E(K_v)$ injects into $\Phi_v(k_v)$.
\end{proof}
\begin{lemma}\label{AppendixLemmaKodairaSymbols}
	Assuming Context \ref{ContextArbitraryGenus}, the reduction types of $E/K$ and of $E_f/K$ are related as follows.
	\begin{enumerate}
		\item[(a)] If $v \in \mid C\mid _{\text{unr}}$,
		\begin{enumerate}
			\item[(a1)] then $v \in \mid U\mid  \cap \mid U_f\mid $ or $v \in M \cap M_f$ or $v \in A \cap A_f$.
			\item[(a2)] Moreover, if $v \in M$, then $v \in M^{\text{sp}} \cap M_f^{\text{sp}}$ or $M^{\text{ns}} \cap M_f^{\text{ns}}$ if and only if $v \in M_{\text{split}}$.
		\end{enumerate}  
		\item[(b)] If $v \in \mid C\mid _{\text{ram}}$, then
		\vspace{0.4cm}	
		\begin{table}[h]
			\centering
			\begin{tabular}{ |c|c|c|c|c|c| }
				\hline
				$E/K_v$ & $\mid U\mid _{\text{ram}}$ & $M_{\text{ram}}$ & $A_\text{ram}^{\text{gd}}$ & $A_\text{ram}^m$ & $A_\text{ram}^{o}$\\
				\hline
				$E_f/K_v$ & $A_{\text{ram}}^\text{gd}$  & $A_{f,\text{ram}}^m$ & $\mid U_f\mid _{\text{ram}}$ & $M_{f,\text{ram}}$ & $A_{f,\text{ram}}^{o}$\\
				\hline
			\end{tabular}.
			\caption{Reduction Types After Ramified Twist}
			\label{tb22RamifiedTwist}
		\end{table}
	
	\end{enumerate}
\end{lemma}
\begin{proof}
	Suppose that the elliptic curve $E/K_v$ is given by a minimal Weierstrass equation
	\begin{equation}\label{InitialWeierstrassEquation}
		y^2 = x^3 + a_v x + b_v
	\end{equation}
	and let $\Delta_v$ be the discriminant associated to this equation. Then, $E_{f_v}/K_v$ is given by the Weierstrass equation
	\begin{equation}\label{WeierstrassEquationTwist}
		y^2 = x^3 + f_v^2 a_v x + f_v^3 b_v,
	\end{equation}
	and the discriminant of this equation is $f_v^6 \Delta_v$ (see Lemma \ref{DiscriminantAndJInvariantQuadraticTwist}). However, Equation (\ref{WeierstrassEquationTwist}) might not be minimal. We will reduce to the case where it is minimal.
	
	There exists an integer $r \in \Z$ such that $f_v \pi_v^{-2r} \in \calO_v$. Under the change of variables $(x,y) \mapsto (x/\pi_v^{2r}, y/\pi_v^{3r})$ the Weierstrass equation (\ref{WeierstrassEquationTwist}) becomes
	\begin{equation}\label{WeierstrassEquationTwistModified}
		y^2 = x^3 + f_v^2\pi_v^{-4r}a_v x + f_v^3\pi_v^{-6r}b_v
	\end{equation}
	and has all its coefficients in $\calO_v$ and therefore the discriminant of Equation (\ref{WeierstrassEquationTwistModified}), $f_v^6\pi_v^{-12r}\Delta_v$, also belongs to $\calO_v$. Since $\ord_v$ is a discrete valuation, we can choose $r \in \Z$ such that 
	$$
	\ord_v(f_v^6\pi_v^{-12r}\Delta_v) = 6\left( \ord_v(f_v) - 2r\right) + \ord_v(\Delta)
	$$
	is minimized under the constraint that the coefficients of the Weierstrass equation (\ref{WeierstrassEquationTwistModified}) belong to $\calO_v$, i.e, under the constraint that the integers
	$$
	2\left( \ord_v(f_v) - 2r\right) + \ord_v(a_v) \text{ and } 3\left( \ord_v(f_v) - 2r\right) + \ord_v(b_v)
	$$
	are nonnegative. For any such $r$, the Weierstrass equation (\ref{WeierstrassEquationTwistModified}) is minimal.
	
	Because twisting by an element in the square class of $f_v$ yields an elliptic curve which is isomorphic to $E_{f_v}/K_v$ (see Subsection \ref{QuadraticTwists}), we can assume without loss of generality that the Weierstrass equation (\ref{WeierstrassEquationTwist}) is minimal and that either $\ord_v(f_v) = 0$ (and so $f_v$ is a unit) or that $\ord_v(f_v) \in \{ \pm 1 \}$ with $f_v = \pi_v^{\pm 1}$. Recall from Context \ref{ContextArbitraryGenus} that $K_{f,w} \simeq K_v(\theta_v)$. The integral closure $\calO_w$ of $\calO_v$ inside $K_{f,w}$ is also the valuation ring of the place $w$ lying over $v$ and we let $k_w$ be the residue field of $w$.
	
	Suppose until further notice that $f_v$ is a unit. Then, $\ord_v(f_v) = 0$ and so the extension $K_{f,w}/K_v$ is unramified. Therefore, the canonical surjective group morphism $\Gal(K_{f,w}/K_v) \to \Gal(k_w/k_v)$ is an isomorphism (see \cite[p.172, II, (9.9) Proposition]{neukirch_1999}).
	
	From this observation one sees that $f_v$ is a square in $K_v$ if and only if $v$ splits in $K_{f,w}$ or equivalently that $f_v \bmod \pi_v$ is a square in $k_v$. As $f_v$ is a unit, it follows that
	$$
	\ord_v(f_v^6 \Delta_v) = \ord_v(\Delta_v), \ord_v(f_v^2a_v) = \ord_v(a_v) \text{ and }\ord_v(f_v^3b_v) = \ord_v(b_v).
	$$
	Then $E/K_v$ and $E_f/K_v$ have the same reduction type (good, multiplicative, additive) over $v$ (see \cite[p.196, VII.5.1]{silverman_2009}). This proves Lemma \ref{LemmaKodairaSymbols} (a1). 
	
	Suppose that $v \in M_{\text{unr}}$. By \cite[p.364, Lemma 2.1]{baig_hall_2012}, $E_f/K_v$ has split multiplicative reduction over $v$ if and only if the image $f_v^36b_v \bmod \pi_v$ of $f_v^3b_v$ in $k_v$ is a square. This is equivalent to require that $f_v6b_v \bmod \pi_v$ is a square in $k_v$. But from \cite[p.364, Lemma 2.1]{baig_hall_2012} we also know that $E/K_v$ has split multiplicative reduction over $v$ if and only if $6b_v \bmod \pi_v$ is a square in $k_v$. Hence, $E/K_v$ and $E_f/K_v$ have the same splitting behaviour (both split multiplicative or both nonsplit multiplicative) over $v$ if and only if $f_v \bmod \pi_v$ is a square in $k_v$. This completes the proof of Lemma \ref{LemmaKodairaSymbols} (a2).
	
	We now assume that $f_v \in \{ \pi_v^{\pm 1} \}$. In this situation, $v \in \mid C\mid _{\text{ram}}$. Remember that the Weierstrass equations (\ref{InitialWeierstrassEquation}) and (\ref{WeierstrassEquationTwist}) are minimal.  
	
	Suppose until further notice that $\ord_v(f_v) = 1$. Since Equation (\ref{WeierstrassEquationTwist}) is minimal, then \cite[p.186, VII.1.1]{silverman_2009} implies (not necessarily mutually exclusively) that 
	\begin{enumerate}
		\item[(b1)] $\ord_v\left( \Delta_{v,f}\right) = 6 + \ord_v(\Delta_v) < 12$, which implies that $0 \leq \ord_v(\Delta_v) < 6$ (where the left inequality holds by the minimality of Equation (\ref{InitialWeierstrassEquation})), or
		\item[(b2)] $\ord_v\left( f_v^2a_v \right) = 2 + \ord_v(a_v) < 4$, which implies that $0 \leq \ord_v(a_v) < 2$ (where again the left inequality holds by the minimality of Equation (\ref{InitialWeierstrassEquation})).
	\end{enumerate}
	Suppose until further notice that (b1) holds. 
	\begin{enumerate}
		\item[(b1.1)] If $\ord_v(\Delta_v) = 0$, then $E/K$ has good reduction over $v$ (see \cite[p.196, VII.5.1(a)]{silverman_2009}). In this case, $\ord_v(\Delta_{v,f}) = 6 > 0$ and so $E/K_f$ has bad reduction over $v$. Since $j(E_f) = j(E)$, then
		$$
		\ord_v\left(j(E_f)\right) = 3\ord_v(a_v) - \ord_v(\Delta_v) = 3\ord_v(a_v) \geq 0
		$$
		by the minimality of Equation (\ref{InitialWeierstrassEquation}). Therefore, $j(E_f) \in \calO_v$ and so $E_f/K_f$ has (additive) potential good reduction over $v$ (see \cite[p.197, VII.5.5]{silverman_2009}).
		\item[(b1.2)] If $\ord_v(\Delta) > 0$, then $E/K$ has bad reduction over $v$.
		\begin{enumerate}
			\item[(b1.2.1)] If $\ord_v(a_v) = 0$, then $E/K$ has multiplicative reduction over $v$ by \cite[p.196, VII.5.1(b)]{silverman_2009} and in this case $\ord_v(f_v^2a_v) = 2 > 0$. This implies that $E_f/K$ has additive reduction over $v$ (see \cite[p.196, VII.5.1(c)]{silverman_2009}).
			\item[(b1.2.2)] Suppose that $\ord_v(a_v) > 0$. Then $E/K$ has additive reduction over $v$ (see \cite[p.196, VII.5.1(c)]{silverman_2009}). Since $\ord_v(f_v^2a_v) = 2 + \ord_v(a_v) > 0$, then $E_f/K$ also has additive reduction over $v$ (see \cite[p.196, VII.5.1(c)]{silverman_2009}).
		\end{enumerate}
	\end{enumerate}
	Now suppose that (b2) holds. Our case by case analysis of (b1) covers the following situations: When $\ord_v(\Delta_v) = 0$, we conclude as in (b1.1) for $\ord_v(a_v) \in \{0,1\}$. When $\ord_v(\Delta_v) > 0$ and $\ord_v(a_v) = 0$, then we conclude as in (b1.2.1). When $0 < \ord_v(\Delta_v) < 6$ and $\ord_v(a_v) = 1$, we conclude as in (b1.2.2). We are left to consider the situation where $\ord_v(\Delta_v) \geq 6$ and $\ord_v(a_v) = 1$. In this situation $E/K$ has additive reduction over $v$ and so does $E_f/K$.
	
	Finally, assume that $\ord_v(f_v) = - 1$. Then, $\ord_v(\Delta_{v,f}) \geq 0$ (resp. $\ord_v(f_v^2a_v)) \geq 0$) implies that $\ord_v(\Delta_v) \geq 6$ (resp. $\ord_v(a_v) \geq 2$). Since Equation (\ref{InitialWeierstrassEquation}) is minimal, then $E/K$ has additive reduction over $v$. Since Equation (\ref{WeierstrassEquationTwist}) is minimal, then \cite[p.186, VII.1.1]{silverman_2009} implies that we have the following (not necessarily mutually exclusively) two possibilities.
	\begin{enumerate}
		\item[(b3)] If $\ord_v(\Delta_{v,f}) = - 6 + \ord_v(\Delta_v) < 12$, then $6 \leq \ord_v(\Delta_v) < 18$.
		\item[(b4)] If $\ord_v(f_v^2 a_v) = -2 + \ord_v(a_v) < 4$, then $2 \leq \ord_v(a_v) < 6$.
	\end{enumerate}
	Suppose that (b3) holds.
	\begin{enumerate}
		\item[(b3.1)] If $\ord_v(\Delta_{v,f}) = 0$, then $\ord_v(\Delta_v) = 6$. Since $E_f/K$ has good reduction over $v$, then $j(E_f) \in \calO_v$ by \cite[p.197, VII.5.5]{silverman_2009}. We also have $j(E_f) = j(E)$. Therefore, 
		$$
		\ord_v(j(E)) = 3\ord_v(a_v) - 6 \geq 0
		$$
		and so $\ord_v(a_v) \geq 2$. Thus, by \cite[p.197, VII.5.5]{silverman_2009} the elliptic curve $E/K$ has (additive) potential good reduction over $v$.
		\item[(b3.2)] If $\ord_v(\Delta_{v,f}) > 0$, then $E_f/K$  has bad reduction over $v$. Now, $\ord_v(a_v) \geq 2$.
		\begin{enumerate}
			\item[(b3.2.1)] If $\ord_v(a_v) = 2$, then $\ord_v(f_v^2a_v) = 0$ and so $E_f/K$ has multiplicative reduction over $v$ by \cite[p.196, VII.5.1(c)]{silverman_2009} and $E/K$ has additive reduction. More precisely, if $\ord_v(j(E_f)) = -\ord_v(\Delta_{v,f}) = -n$ for some positive integer $n$, so that $E_f/K$ has Kodaira symbol $I_n$ over $v$, then $\ord_v(j(E)) = -n$ as well and $E/K$ has potential multiplicative reduction over $v$ with Kodaira symbol $I_n^*$.
			\item[(b3.2.2)] If $\ord_v(a_v) > 2$, then $E_f/K$ has additive reduction over $v$. 
		\end{enumerate}
	\end{enumerate}
	Finally, suppose that (b4) holds. If $\ord_v(\Delta_{v,f}) = 0$, then we conclude as in (b1.1). Suppose that $\ord_v(\Delta_{v,f}) > 0$. If $\ord_v(a_v) = 2$, then we conclude as in (b3.2.1), while if $\ord_v(a_v) > 2$, then we conclude as in (b3.2.2).
	
	This ends our case by case analysis and our proof of Lemma \ref{AppendixLemmaKodairaSymbols}.
\end{proof}
\begin{lemma}\label{LemmaLegendreSymbolAppendix}
	For each  $v \in \mid U \mid_\unr = \mid U \mid \cap \mid U_f\mid $, we have $a_v + a_{f,v} \equiv 0 \bmod 2$, where $a_v$ and $a_{f,v}$ are defined in Equation (\ref{TraceOfFrobenius}).
\end{lemma}
\begin{proof}
	For each $v \in \mid U\mid _\unr$, let $\chi_v : k_v^\times \to \{ \pm 1 \}$ be the unique character of order $2$ defined by $\chi_v(r) = 1$ if and only if $r$ is a square in $k_v$, and $\chi_v(r) = -1$, otherwise, and extend $\chi_v$ to $k_v$ by setting $\chi_v(0) := 0$. Moreover, choose a Weierstrass equation $y^2 = \alpha_v(x)$ for the elliptic curve $E_v/k_v$. A corresponding Weierstrass equation for $E_{f,v}/k_v$ is then given by $y^2 = f_v^{-1}\alpha_v(x)$. Here, $f_v$ is identified with the unit $f_v$, which is the image of $f$ in the local ring $\calO_v$ (see the proof of Lemma \ref{AppendixLemmaKodairaSymbols}). By \cite[p.139, V.1.3]{silverman_2009}, we have
	$$
	\#E_v(k_v) = 1 + q^{d_v} + \sum_{x \in k_v} \chi_v(\alpha_v(x)) \text{ and } \#E_{f,v}(k_v) = 1 + q^{d_v} + \chi_v(f_v^{-1})\sum_{x \in k_v} \chi_v(\alpha_v(x)).
	$$
	Since $q \equiv 1 \bmod 2$, then
	$$
	\#E_v(k_v) + \#E_{f,v}(k_v) = 2(1+q^{d_v}) + (1 + \chi_v(f_v^{-1}))\sum_{x \in k_v} \chi_v(\alpha_v(x)) \equiv 0 \bmod 2.
	$$
	In particular, we have
	$$
	a_v + a_{f,v} = 2(1+q^{d_v}) - \left( \#E_v(k_v) + \#E_{f,v}(k_v) \right) \equiv 0 \bmod 2.
	$$
\end{proof}
\section*{Acknowledgments}
This paper comes from Chapter 3 of the PhD thesis of the author, while at the University of Western Ontario. We acknowledge the support of the Natural Sciences and Engineering Research Council of Canada (NSERC), [PGSD2-535451-2019] during part of this period. The paper was completed at the University of Lethbridge where the author received financial support from the Pacific Institute for the Mathematical Sciences (PIMS) while he was a PIMS postdoctoral fellow. The author sincerely thanks his PhD supervisor Chris Hall for the many helpful discussions. The author also thanks Amir Akbary and Andrew Fiori for the useful feedback on preliminary versions of this paper.

\end{document}